\documentclass[a4paper,11pt]{amsart}
\usepackage[top=2.2cm,left=2.8cm,right=2.8cm,bottom=2.8cm]{geometry}

\usepackage[foot]{amsaddr}

\usepackage{hyperref}
\hypersetup{
    colorlinks,
    citecolor=green,
    filecolor=black,
    linkcolor=blue,
    urlcolor=blue
}
\hypersetup{linktocpage}

\usepackage{hhline}

\usepackage{subfloat}
\usepackage{mathtools}
\usepackage{cite}
\usepackage{graphicx}
 
 \usepackage{cleveref}
\usepackage{amsbsy}
\usepackage{latexsym}
\usepackage{amsfonts}
\usepackage{amssymb}
\usepackage{upgreek}
\usepackage[usenames]{color}
\usepackage{amsmath,amsthm}
\usepackage{enumerate}
\usepackage{stmaryrd}

\usepackage{multirow}

\usepackage{pgfplots} 
\usepackage{tikz}
\usepackage{tikzscale}
\usetikzlibrary{shapes,backgrounds,calc}

\newtheorem{prop}{Proposition}[section]

\newtheorem{rem}[prop]{Remark}

\usepackage[utf8]{inputenc}
\usepackage{dsfont}
\usepackage{color}
\usepackage{graphicx}
\usepackage{mdwlist}
\usepackage{enumerate}
\usepackage{setspace}
\usepackage{multicol}
\usepackage{multirow}
\usepackage{cite}
\usepackage{mathtools}

\usepackage[section]{algorithm}
\usepackage{algorithmicx}
\usepackage{algpseudocode}
\algrenewcommand\algorithmicrequire{\makebox[32pt][l]{\textrm{input}}}
\algrenewcommand\algorithmicensure{\makebox[32pt][l]{\textrm{output}}}
\algrenewcommand\algorithmicfunction{\textrm{function}}
\algrenewcommand\algorithmicwhile{\textrm{while}}
\algrenewcommand\algorithmicdo{}
\algrenewcommand\algorithmicend{\textrm{end}}
\algrenewcommand\algorithmicforall{\textrm{for all}}
\algrenewcommand\algorithmicfor{\textrm{for}}
\algrenewcommand\algorithmicrepeat{\textrm{repeat}}
\algrenewcommand\algorithmicuntil{\textrm{until}}

\DeclareFontFamily{U}{mathx}{\hyphenchar\font45}
\DeclareFontShape{U}{mathx}{m}{n}{
      <5> <6> <7> <8> <9> <10>
      <10.95> <12> <14.4> <17.28> <20.74> <24.88>
      mathx10
      }{}
\DeclareSymbolFont{mathx}{U}{mathx}{m}{n}
\DeclareMathSymbol{\bigtimes}{1}{mathx}{"91}

\newtheorem{theorem}{Theorem}[section]
\newtheorem{lemma}[theorem]{Lemma}

\newtheorem{remark}[theorem]{Remark}

\theoremstyle{definition}

\numberwithin{equation}{section}

\newcommand{\N}{\mathds{N}}
\newcommand{\R}{\mathds{R}}
\newcommand{\Z}{\mathds{Z}}
\newcommand{\C}{\mathds{C}}
\newcommand{\Q}{\mathds{Q}}

\DeclareMathOperator*{\argmin}{arg\,min}

\newcommand{\D}{{\rm D}}

\newcommand{\spl}[1]{\ell_{#1}}

\providecommand{\norm}[1]{\lVert#1\rVert}

\def\bu{{\bf u}}

\def\bV{{\bf V}}
\def\bW{{\bf W}}
\newcommand{\bv}{{\bf v}}
\newcommand{\bw}{{\bf w}}
\def\e2{\spl{2}(\nabla^d)}
\newcommand{\cA}{\mathcal{A}}
\def\ga{\gamma}

\newcommand{\bM}{\mathbf{M}}
\newcommand{\bT}{\mathbf{T}}
\newcommand{\bbf}{\mathbf{f}}
\newcommand{\cL}{\mathcal{L}}
\newcommand{\cD}{\mathcal{D}}

\newcommand{\cN}{\mathcal {N}}

\newcommand{\cP}{\mathcal{P}}
\renewcommand{\e}{\varepsilon}
\newcommand{\cU}{{U}}

\newcommand{\bG}{\mathbf{G}}
\newcommand{\T}{\mathbb{T}}

\newcommand{\osc}{{\rm vis}}
\begin{document}
\title{Convergence of Adaptive Least Squares Finite Element Methods}
\thanks{This work has been supported in
 part by NSF Grants  DMS-2513112, DMS-2245097, DMS-2038080,   and by the German Research Foundation grant CRC 1481.
% by  the Office of
%Naval Research Contract N0014-91-J1343, the Army Research Office contract Nr.
%DAAD 19-02-1-0028, the TMR network
%``Wavelets in Numerical Simulation", by the SFB 401, funded by the
%German Research Foundation and by the Alexander von Humboldt- Foundation}
}
\author{Wolfgang Dahmen} %

\address{University of South Carolina, 1523 Greene Street, 29208 Columbia, South
  Carolina}
\email{wolfgang.anton.dahmen@googlemail.com}

\renewcommand{\arraystretch}{1.5}
\hbadness=10000
\vbadness=10000
\newtheorem{propo}{Proposition}[section]
\newtheorem{coro}[prop]{Corollary}
\newtheorem{rema}[prop]{Remark}
\def\ds{\displaystyle}
\def\ev#1{\vec{#1}}     % arrow over letter
\newcommand{\lt}{\ell_{2}(\cJ)}
\def\Supp#1{{\rm supp\,}{#1}}
\def\e{\varepsilon}
\def\nl{\newline}
\def\T{{\relax\ifmmode I\!\!\hspace{-1pt}T\else$I\!\!\hspace{-1pt}T$\fi}}
\def\N{\mathbb{N}}
\def\R{\mathbb{R}}
\def\Z{{\relax\ifmmode Z\!\!\!Z\else$Z\!\!\!Z$\fi}}
\def\Zd{\Z^d}
\def\Q{{\relax\ifmmode I\!\!\!\!Q\else$I\!\!\!\!Q$\fi}}
\def\C{{\relax\ifmmode I\!\!\!\!C\else$I\!\!\!\!C$\fi}}
\def\Rd{\R^d}
\def\gsim{\mathrel{\raisebox{-4pt}{$\stackrel{\textstyle>}{\sim}$}}}
\def\sime{\raisebox{0ex}{$~\stackrel{\textstyle\sim}{=}~$}}
\def\lsim{\raisebox{-1ex}{$~\stackrel{\textstyle<}{\sim}~$}}
\def\div{\mbox{div}}
\def\gr{\nabla}
\def\curl{\mbox{ curl }}
\def\M{M}  \def\NN{N}                  % parameter M, N in [DKU]
\def\L{{\ell}}               % parameter L = l in [DKU] (support)
\def\Le{{\ell_1}}            % support of xi: [\ell_1,\ell_2]
\def\Lz{{\ell_2}}
\def\Let{{\tilde\ell_1}}     % support of xi\tilde:[\t\ell_1,\t\ell_2]
\def\Lzt{{\tilde\ell_2}}
\def\Ltw{\ell_\tau^w(\cJ)}
\def\tLtw{{}_t\Asl}
\def\t#1{\tilde{#1}}
\def\la{\lambda}
\def\La{\Lambda}
\def\ga{\gamma}
\def\Ga{\eta}
\def\al{\alpha}
\def\cO{\mathcal{O}}
\def\cA{\mathcal{A}}
\def\cC{\mathcal{C}}
\def\cF{\mathcal{F}}
\def\bu{{\bf u}}
\def\bbf{{\bf f}}
\def\bz{{\bf z}}
\def\bZ{{\bf Z}}
\def\bI{{\bf I}}
\def\bN{{\bf N}}
\def\cE{\mathcal{E}}
\def\cD{\mathcal{D}}
\def\cI{\mathcal{I}}
\def\cJ{\mathcal{J}}
\def\cM{\mathcal{M}}
\def\cN{\mathcal{N}}
\def\cT{\mathcal{T}}
\def\cU{\mathcal{U}}
\def\cV{\mathcal{V}}
\def\cL{\mathcal{L}}
\def\cB{\mathcal{B}}
\def\cG{\mathcal{G}}
\def\cK{\mathcal{K}}
\def\cS{\mathcal{S}}
\def\cP{\mathcal{P}}
\def\cQ{\mathcal{Q}}
\def\cR{\mathcal{R}}
\def\cU{\mathcal{U}}
\def\bL{{\bf L}}
\def\bK{{\bf K}}
\def\bC{{\bf C}}
\def\X{X\in\{L,R\}}
\def\ph{{\varphi}}
\def\D{{\Delta}}
\def\H{\mathcal{H}}
\def\bM{{\bf M}}
\def\bQ{{\bf Q}}
\def\bG{{\bf G}}
\def\bP{{\bf P}}
\def\bW{{\bf W}}
\def\bT{{\bf T}}
\def\bV{{\bf V}}
\def\bv{{\bf v}}
\def\bz{{\bf z}}
\def\bw{{\bf w}}
\def\bQ{{\bf Q}}
\def\t{\tilde}
\def\lll{\langle}
\def\rr{\rangle}
\def\dJ{\nabla}
\newcommand{\ba}{{\bf a}}
\newcommand{\bb}{{\bf b}}
\newcommand{\bc}{{\bf c}}
\newcommand{\bd}{{\bf d}}
\newcommand{\bs}{{\bf s}}
\newcommand{\bff}{{\bf f}}
\newcommand{\bg}{{\bf g}}
\newcommand{\bq}{{\bf q}}
\newcommand{\br}{{\bf r}}
\newcommand{\beze}{{\bf 0}}
\newcommand{\PP}{\mathbb{P}}
\newcommand{\gkl}{\Gamma_{kl}}
\newcommand{\Hoo}{H^{1/2}_{00}}
\newcommand{\Hoogkl}{\Hoo(\gkl)}
\newcommand{\Asl}{{\cA^s}}
\newcommand{\AsH}{{\cA^s(\cH)}}
\newcommand{\Aslt}{{\cA^s_{\rm tree}}}
\newcommand{\AsHt}{{\cA^s_{\rm tree}(\cH)}}
\newcommand{\beqn}{\begin{equation}}
\newcommand{\eeqn}{\end{equation}}
\newcommand{\be}{\begin{equation}}
\newcommand{\ee}{\end{equation}}
\def\beginproof{\noindent{\bf Proof:}~ }
\def\endproof{\hfill\rule{1.5mm}{1.5mm}\\[2mm]}

% -----------------------------
% Neue Definitionen (Th. Lange)
% -----------------------------
\newenvironment{Proof}{\noindent{\bf Proof:}\quad}{\endproof}

\renewcommand{\theequation}{\thesection.\arabic{equation}}
\renewcommand{\thefigure}{\thesection.\arabic{figure}}

\makeatletter
\@addtoreset{equation}{section}
\makeatother

\newcommand\clos{\mathop{\rm clos}\nolimits}
\newcommand\trunc{\mathop{\rm trunc}\nolimits}
\renewcommand\d{d}

\newcommand\diam{\mathop{\rm diam}}
\newcommand\cond{\mathop{\rm cond}\nolimits}
\newcommand\eref[1]{(\ref{#1})}
\newcommand\Hnorm[1]{\norm{#1}_{H^s([0,1])}}
\def\int{\intop\limits}
\renewcommand\labelenumi{(\roman{enumi})}
\newcommand\lnorm[1]{\norm{#1}_{\ell_2(\Z)}}
\newcommand\Lnorm[1]{\norm{#1}_{L_2([0,1])}}
\newcommand\LR{{L_2(\R)}}
\newcommand\LRnorm[1]{\norm{#1}_\LR}
\newcommand\Matrix[2]{\hphantom{#1}_#2#1}
\newcommand\ogauss[1]{\left\lceil#1\right\rceil}
\newcommand{\QED}{\hfill \raisebox{-2pt}{\rule{5.6pt}{8pt}\rule{4pt}{0pt}}%
  \smallskip\par}
\newcommand\Rscalar[1]{\scalar{#1}_\R}
\newcommand\scalar[1]{\left(#1\right)}
\newcommand\Scalar[1]{\scalar{#1}_{[0,1]}}
\newcommand\Span{\mathop{\rm span}}
\newcommand\ugauss[1]{\left\lfloor#1\right\rfloor}
\newcommand\with{\, : \,}

\newcommand\bnull{{\bf 0}}
\newcommand\bA{{\bf A}}
\newcommand\bB{{\bf B}}
\newcommand\bR{{\bf R}}
\newcommand\bD{{\bf D}}
\newcommand\bE{{\bf E}}
\newcommand\bF{{\bf F}}
\newcommand\bH{{\bf H}}
\newcommand\cH{\mathcal{H}}
\newcommand\sk{{sk}}
\newcommand{\XX}{\mathbb{X}}
\newcommand{\id}{{\rm id}}

%-----------------------
% Neues von T. Bronger
%-----------------------
\newif\ifNZB
\newcommand\NZB[1]{\ifNZB \marginpar{\raggedright \scriptsize NZB:\\#1}
 \else \fi}
\newcommand{\FText}[1]{\mbox{#1}}
\makeatletter
\newcommand{\tr}{{\mathop{\operator@font T}\nolimits}}
\makeatother
\newcommand{\UArrow}[4]{
 \begin{array}{ll}
  #1&\stackrel{#2}\longrightarrow\\
  #3&\;\raisebox{1ex}{$\nearrow$}\mkern-14mu_{#4}
 \end{array}}
\newcommand{\DArrow}[4]{
 \begin{array}{ll}
  \stackrel{#1}\longrightarrow&#2\\
  \mkern-10mu_{#3}\mkern-22mu\raisebox{1ex}{$\searrow$}&#4
 \end{array}}
\newcommand{\fig}[3]{\par\begin{figure}[ht]
  \centerline{\epsfbox{#1.eps}}\caption{#3}\label{fig#2}\end{figure}}

\newcommand{\ubk}{u^{\rm bk}}
\newcommand{\utr}{u^{\rm true}}
\newcommand{\Atr}{A^{\rm true}}
%-----------------------
%%%%%%%%%%%%%%%%%%%%%%%%%%%%%von mir
\newcommand{\dI}{\Delta}
\newcommand{\cY}{\mathcal{Y}}
\newcommand{\U}{\mathbb{U}}
\newcommand{\V}{\mathbb{V}}
\newcommand{\W}{{\mathbb{W}^\ell}}
\newcommand{\Wi}{{\mathbb{W}^\infty}}
\newcommand{\pP}{\mathbb{P}}
\newcommand{\Um}{\mathbb{U}_{\rm mixed}}
\newcommand{\Vm}{\mathbb{V}_{\rm mixed}}
\newcommand{\oVm}{\tilde{\mathbb{V}}_{\rm mixed}}
\newcommand{\Y}{\mathbb{Y}}

\newcommand{\E}{\mathbb{E}}
\newcommand{\Lip}{{\rm Lip}}
\newcommand{\cX}{\mathcal{X}}

\def\argmin{\mathop{\rm argmin}}
\def\nargmin{\mathop{\rm near-argmin}}

\def\argmax{\mathop{\rm argmax}}

 \newcommand{\tripnorm}[1]{|\!|\!| #1|\!|\!|}
 \newcommand{\w}{\mathsf{w}}
 
 \newcommand{\dN}{\mathcal{N}\!\mathcal{N}}
\newcommand{\cZ}{\mathcal{Z}}
\renewcommand{\tr}{{\rm train}}
\newcommand{\hx}{{\hat{x}}}
\renewcommand{\H}{\mathbb{H}}
\newcommand{\s}{\mathsf{s}}
\newcommand{\utau}{\underline{\tau}}
\newcommand{\otau}{\overline{\tau}}
%%%%%%%%%%%%%%%%%%%%%%%
\newcommand{\Ud}{{\U^\delta}}
\newcommand{\Vd}{{\V^\delta}}
\newcommand{\hVd}{{\hat{\V}^\delta}}

\newcommand{\mD}{\mathfrak{D}}
\newcommand{\wt}{\widetilde}
\renewcommand{\X}{\XX}

\newcommand{\mA}{\mathfrak{A}}
\newcommand{\mI}{\mathfrak{I}}
\newcommand{\mS}{\mathfrak{S}}
\newcommand{\mR}{\mathfrak{R}}
\newcommand{\mP}{\mathfrak{P}}
\newcommand{\md}{\mathfrak{d}}

\newcommand{\wtF}{{\cF}}
\newcommand{\wtX}{{\X}}
\newcommand{\wtY}{{\Y}}
\newcommand{\wtB}{{\cB}}
\newcommand{\wh}{\widehat}

\newcommand{\pp}{y} %{\textsc{p}}
\newcommand{\pdom}{\cY}%{\mathcal{P}}
\newcommand{\ccY}{\pdom} 

\newcommand{\uN}{{u^{\mbox{\tiny{$\mathcal{N}\!$}\! }}}}
\newcommand{\uNe}{{u_\e^{\mbox{\tiny{$\mathcal{N}\!$}\! }}}}
\newcommand{\uNi}{{u_i^{\mbox{\tiny{$\mathcal{N}\!$}\! }}}}

\renewcommand{\mP}{\Theta}

\newcommand{\loss}{\cL}
\newcommand{\lossemp}{\loss_{\rm emp}}
\newcommand{\hu}{\widehat{u}}
%%%%%%%%%%%%%%%%%%%
%%%%%%%%%%%%%%%%%%%%%%%%%%%%%%%%%%%%%%%%%%%%%%
\newcommand{\Orc}{{\Omega^\circ_{\rm ref}}}
\newcommand{\Or}{{\Omega_{\rm ref}}}
\newcommand{\Ori}{{\Omega_{{\rm ref},i}}}
\newcommand{\Ur}{\U_{\rm ref}}

\newcommand{\cKo}{{\cK_{\rm osc}}}
\newcommand{\cKr}{{\cK_{\rm res}}}
\newcommand{\cTr}{{\cT_{\rm ref}}}
\newcommand{\Pref}{P_{\U_\cTr}}
\newcommand{\wT}{w_\cT}

\newcommand{\Om}{\Omega}
\newcommand{\UTr}{{\U_{\cTr}^\infty}}

\newcommand{\wTc}{u^{\rm ref}_{\cT}}
\newcommand{\wi}{{w^\circ}}

\newcommand{\inv}{{\rm inv}}

\newcommand{\wc}{w^\circ}
%%%%%%%%%%%%%%%%%%%%%%%%%%%%%%%%

%
%\date{}
\maketitle
\begin{abstract}
We consider a general class of operator
equations $\cB u=f$ where $\cB:\U\to \mathbb L$ is a linear isomorphis from a Hilbert space $\U$
onto an ``$L_2$-type'' space $\mathbb L = L_2(\Omega;\R^q)$. This covers weak first order formulations of stationary elliptic boundary value problems, the Stokes system or systems for linear elasticity
as well as space-time formulations of heat or wave equations. Since an $L_2$-residual is uniformly
equivalent to the corresponding error in the trial space norm reliable and efficient computable 
a posteriori error bounds are for free. In contrast to a matured theory of classical adaptive finite element methods, estblishing linear convergence - and eventually optimal complexity - in a least squares
setting seems to have remained open. To a large extent this is due to the fact that least squares 
error indicators do not exhibit any explicit dependence on local
mesh-sizes. In this paper we propose a rather different approach to establish under which circumstances
standard D\"orfler-marking mesh refinement gives rise to linear convergence  and rate optimality of adaptive least squares
finite element methods as well. A central role is played by a new notion of {\em data-visibility}
which is a structural property of the data rather than a smoothness condition.
\end{abstract}

\vspace*{2mm}

 %%%%%%%%%%%%%%%%%%%%%%%%%%%%%
\section{Introduction}
%%%%%%%%%%%%%%%%%%%%%%%%%%%%%%%%
\subsection{Least Squares Formulations}
%%%%%%%%%%%%%%%%%%%%%%%
{\em Well-posedness} of a (formal) linear operator equation 
\be
\label{opeq}
\cB u =f
\ee 
 is inevitably tied to the specification
of a (linear) trial space $\U$ - often a Hilbert space - that is supposed to accommodate a  unique solution $u^*\in\U$, and a test space $\V$, whose dual $\V'$ then determines the admissible range for the source $f$. Existence of a solution for all $f\in \V'$, its uniqueness and continuous
 dependence  on the data $f$ then means that $\cB:\U\to \V'$ is a norm-isomorphism. 
 This   is well known to be equivalent to the fact that the bilinear form $b(u,v):= (\cB u)(v)$,
 that takes $\U\times \V$ into $\R$, satisfies the conditions of the Babu\v{s}ka-Ne\v{c}as Theorem (see Section \ref{ssec:FOSLS} for more details). This, in turn means that there exist constants $0<c_b\le C_b<\infty$ such that
\be
\label{err-res0}
c_b\|u- u^* \|_\U\le \|f- \cB u\|_{\V'}\le C_b\|u- u^* \|_\U\quad \forall\, u\in\U.
\ee
Thus, on the one hand,   the residual in $\V'$ is a reliable and efficient error estimator,
and on the other hand,   the exact solution $u^*$ as well as an approximation to $u^*$ from some subspace
$\U_n$, can be found as the minimizer of the respective (squared) residual. 

When \eref{opeq} is a homogeneous second order elliptic Dirichlet problem,  one finds $\U=H^1_0(\Omega)=\V$
and $\V'=(H^1_0(\Omega))'=: H^{-1}(\Omega)$. This has triggered the development of an art to
derive computable equivalent expressions for $\|f- \cB u\|_{H^{-1}(\Omega)} $ in terms of  
  so-called   dual residual a posteriori error indicators, \cite{V}  or, from a different angle,
approaches based on the hyper-circle principle, see e.g. \cite{BPS}. This   has by now led to a rich theory of {\em Adaptive Finite Element Methods} (\textsc{Afem}) and a mature understanding of the convergence and complexity properties of \textsc{Afem}, see e.g. 
\cite{BDD,axioms,Feischl,KV,MNS,NSV,robfocm} for a rather incomplete list, referring to \cite{BCNVacta}
for a state-of-the-art overview and further references. Specifically, the essential structure of
convergence and complexity proofs have been abstracted in \cite{axioms}, providing a roadmap for treating a
wide spectrum of application scenarios.

While in the above context a central challenge (systematically taken up first in \cite{V}) has been to  evaluate the dual  norm $\|f- \cB u\|_{\V'}$. This issue disappears
 when $\V=\V'$ happens to be self-dual, hence some kind of an $L_2$-space (e.g. $L_2(\Omega;\R^q)$ for some $q\in\N$). In fact, computing a residual in $L_2$ (up to data-oscillation) 
is conceptually much simpler than in a non-trivial dual space. In other words, reliable and efficient error estimators are for free,
when $\cB$  happens to be an isomorphism whose range is some $L_2$-space.

Specializing \eref{err-res0}, the solution $u^*$ of \eref{opeq}  is in this case characterized by
\be
\label{min}
u^*=\argmin_{u\in\U}\|f- \cB u\|^2_{L_2},
\ee
coining the term {\em least squares method}.
Regarding the scope of its applicability, the paper \cite{FOSLS} has marked an important breakthrough in proving that a specific weak first order formulation of second order elliptic problems is indeed well-posed and equivalent to \eref{min},
see also \cite{Bochev} and the literature cited there.
 This has meanwhile been extended to other important PDE models, including first order formulations of the system for linear elasticity, of the Stokes system, as well as, more recently, space-time first order formulations of parabolic problems and
 of the acoustic wave equation, \cite{FK,FGK,GS}. 
 
 Such $L_2$-minimum residual methods have by now not only proven meritious in   classical finite elemet settings but more recently  also in the context of operator learning, \cite{BDO,CDG,QCD}. 

It should not surprise that, given the local error indicators $\eta_K(u):=\int_K|f(x)-\cB u(x)|^2dx$ for free, where 
$K$ is an element of a current mesh (assuming a sufficiently accurate evaluations of such integrals for the moment),  corresponding {\em adaptive least squares methods} (\textsc{Alsm})
 have been lately  attracting considerable attention, last but not least, because, empirically, they indeed turn out to perform well, see e.g. \cite{Bringmann} for an overview.
%%%%%%%%%%%%%%%%%%%%%%%%%%%%
\subsection{Why Is there a Challenge for Adaptivity?}
%%%%%%%%%%%%%%%%%%%%%%%%%%%%%%%%%% 
  It may surprise though that the theoretical underpinning of \textsc{Alsm} is by far less mature than for AFEM. Specifically,
 the central issue, namely the question whether  \eref{cycle} entails the reduction of the current error 
 by a fixed factor less than one (which is usually also the basis for proving optimal convergence/complexity rates), appears to be wide open.
The only exception is \cite{CP} where 
for special lowest order discretizations  for the FOSLS formulation of Poisson's equation convergence could
be quantified by emulating conventional techniques for elliptic problems. Unfortunately, this doesn't
seem to lead to any approach that works for more general discretizations and a wider spectrum of problems.   

On a superficial level the main obstruction is simple.
The considerable body of literature and expertise concerning standard \textsc{Afem} for  problems of elliptic type makes essential use of the fact
that corresponding local error indicators involve explicit {\em local mesh-size} factors
when proving that for each adaptive cycle
\be
\label{cycle}%\begin{center}
\textsc{Mark} \quad\to\quad \textsc{Refine} \quad\to\quad \textsc{Solve},
\ee 
the current error is reduced by a {\em fixed factor} $\rho<1$, see e.g. \cite{axioms}. As mentioned earlier, this 
{\em linear} convergence is the indispensible 
first step for eventually proving, in a second step, that the convergence is even optimal
over natural classes of approximation spaces, see e.g. \cite{BDD,BCNVacta,MNS,robfocm}.
Least squares residuals and resulting error indicators {\em do not involve} such mesh-size factors.
%%%%%%%%%%%%%%%%%%%%%%%%%%
\subsection{Contribution}
%%%%%%%%%%%%%%%%%%%%%%%%%%%%%
The central objective of this paper is therefore to shed light on whether nevertheless
suitable versions of \eref{cycle}  give rise to a {\em fixed error reduction per cycle}  also for \textsc{Alsm}.
We do so not for a specific PDE but on an abstract level for \eref{opeq} when it is well-posed for 
a self-dual test space $\V=\V'$ of $L_2$-type.

 To destill the essential mechanisms, we analyze first
the convergence of \eref{cycle} under {\em idealized} assumptions. In essence, this
means that integrals can be computed   and linear systems are solved exactly. This already 
brings out a notion of {\em data visibility} that governs error reduction. 
%perhaps albeit in a slightly different
%way than usually, see \cite{KV}. 
It should be stressed that this new notion, despite some   similarities with what is commonly referred to as {\em data-oscillation}, 
    differs in essential ways.
Among other things, it is not only concerned with the ability to evaluate data exactly but is intimately tied to the current refinement region, made up by the refined cells
 within a given cycle \eref{cycle}. In particular, its role stresses the fact that without any further condition on the right hand side - beyond membership in an $L_2$-type space - convergence performance cannot be quantified. Data-visibility in the present context can be seen as a structural alignment property
 of $f$. Beeing of higher order relative to the residual is not merely a consequence of
smoothness.

These 
findings become most transparent and convenient to handle
when working with a specific norm for the trial space, which is often referred to as {\em optimal trial  norm}.

Once linear error reduction has been established, which is the central issue in this paper, 
following the reasoning in \cite{robfocm}, more practical variants can be handled by perturbation 
arguments, leading then to (optimal) complexity estimates.

In view of the envisaged level of generality, there must remain though some   level of idealization when 
assessing  the {\em complexity} of such methods. We focus here only on reducing  the level of idealization by facilitating the evaluation of residuals via piecewise polynomial approximations to
the load $f$
for the following reasons. First, including this form of {\em perturbing problem data},   is essential to understand error reduction through the previously mentioned notion of data visibility. 
Once, this scenario is under control, further perturbation effects,
like inexact solves or consistency errors due to operator-coefficient approximations, can be handled in essentially the same way as for classical \textsc{Afem}, see e.g. \cite{robfocm}. We therefore
dispense here with repeating this machinery. Therefore, we confine the ``complexity'' discussion
to the growth of partitions which represent  the growth of the number of unknowns in adaptive discretizations. 
Further complexity estimates would require, in the preset generality, just making assumptions
because e.g. solver performance will depend on the nature of a specific problem at hand.
The second reason is that perturbing just $f$ suffices to bring out those parts of the reasoning
that deviate from the complexity arguments, used in classical \textsc{Afem}.

A third reason for discussing at least the present light version of a complexity analysis (rather than
just referring to the many similarities) is the following noteworthy distinction between respective results. While in classical \textsc{Afem} analysis, linear convergence can typically be established 
under no restrictions on the bulk parameter $\theta\in (0,1)$ in D\"orfler marking (other than being less than one),
the complexity analysis requires imposing ``smallness'' conditions on the bulk parameter $\theta$.
This doesn't seem to be necessary  for \textsc{Alsm}.
 %
 %method for the underlying least squares energy.
%
%%%%%%%%%%%%%%%%%%%%%%%%%
\subsection{Layout of the Paper}
%%%%%%%%%%%%%%%%%
 In Section \ref{sec:2} we describe the abstract problem setting
so as to cover the instances mentioned above. We briefly recall as a prototype example
the FOSLS case \cite{FOSLS}. We also  
 formulate the notion of {\em optimal trial norm} in the present context since it will 
 turn out to simplify the convergence analysis. Section \ref{sec:prel} is devoted to collecting the main
 ingredients of an \textsc{Alsm} and to setting corresponding notation. We then proceed formulating
 a  scheme under certain idealized assumptions, focusing on what we call {\em adaptive increment}.
 Its discussion reveals the key obstruction to linear error reduction which leads to the aferomentioned 
  notion of {\em data visibility}.  This is made precise in a first main result whose proof is given in the remainder of Section \ref{sec:prel}. 
  
  In Section \ref{sec:practical} we
 reduce the idealizing assumptions on the right hand side data because this already brings out the 
 essence of  quantifiable convergence under more practical assumptions. In brief,
 quantifying piecewise polynomial approximabiity of the right hand side, allows one to quantify the 
 performnce of the \textsc{Alsm} under canonical assumptions on the trial spaces. These findings are
 collected in the second main result Theorem \ref{thm:pract1}. Finally, in Section \ref{sec:compl}
 we analyze the complexity and convergence under an approximability assumption on the right hand side
 in terms of the growth of adaptively generated partitions and hence the number of degrees of freedom.

\section{Problem Setting}\label{sec:2}
%%%%%%%%%%%%%%%%%%%%%%%%%%%%%%%%%%%%%%%%%%%%%
\renewcommand{\L}{{\mathbb{L}_2}}
%%%%%%%%%%%%%%%%%%%%%%%%%%%
\subsection{FOSLS}\label{ssec:FOSLS}
%%%%%%%%%%%%%%%%%%%%%%%
Let $\U$ be a real Hilbert space with norm $\|\cdot\|$ and $\cB\in \cL(\U,\V)$ 
%is an isomorphism, 
where
$\cL(\U,\V)$ is the space of linear operators from $\U$ to $\V$.
The Babu\v{s}ka-Ne\v{c}as Theorem states that $\cB$ is an {\em isomorphism} (boundedly invertible) if and only if 
the bilinear form $b(u,v):= (\cB u)(v):\U\times\V
\to\R$ is continuous and satisfies (see e.g. \cite{Braess,DGacta}) %an in-sup and surjectivity condition 
\be
\label{BN}
\begin{array}{c}
\displaystyle\inf_{u\in\U}\sup_{v\in\V}\frac{b(u,v)}{\|u\|  \|v\|_\V}\ge c_b  ,\quad 
\sup_{u\in\U}\sup_{v\in\V}\frac{b(u,v)}{\|u\|  \|v\|_\V} \le C_b,
\\[3mm]
\mbox{for $v\in \V$ there exists a $w_v\in \U$ such that $b(w_v,v)\neq 0$}.
\end{array}
\ee
This means that
\be
\label{opbounds}
\|\cB\|_{\cL(\U,\V')}\le C_b,\quad \|\cB^{-1}\|_{\cL(\V',\U)}\le c_b^{-1},
\ee
which in turn says that  for any $f\in\V'$, $\cB u=f$,   in the sense that
\be
\label{weak}
b(u,v)= f(v),\quad v\in \V,
\ee
has a unique solution $u^*$ which is stable,
\be
\label{stab}
C_b^{-1}\|f\|_{\V'}\le \|u^*\|\le c_b^{-1}\|f\|_{\V'}.
\ee
Unless  \eref{weak} arises through an energy minimization, where $\U=\V$ is the natural choice and   coercivity   implies the inf-sup and surjectivity condition, identifying such stable
pairs $\U$, $\V$, may be part of the problem, and one may have to acceppt the fact that \eref{BN} can only hold for $\U\neq \V$.

In fact, in what follows, we focus on this latter scenario when  $\V=\L$ is
an {\em $L_2$-type space}, generally being of the form $L_2(\Omega;\R^q)$ for some 
$q\in\N$ (related to the spatial dimension $d$ of the computational domain $\Omega\subset\R^d$) %.appearing as an outer product of $L_2$-spaces,
hence is {\em self-dual} $\L=\L'$. We postpone describing instances where this is the case
and point out first its relevance.

For convenience we denote in what follows by $(\cdot,\cdot)_\Omega:\L\times \L\to \R$ the 
canonical inner product on $\L$ but prefer to refer to the domain $\Omega$, as a subscript, as this will change at times.
Accordingly, we write
$$
\|g\|_\L = \|g\|_\Om= (g,g)_\Om^{1/2},
$$
because such norms will sometimes be taken with respect to subdomains of $\Om$.

We will make repeated use of the following facts which are in essence known, but stated 
in the following form for the convenience of the reader.
\begin{prop}
\label{prop:equiv}
Let $\U$ be a real Hilbert space of functions defined on a computational domain $\Omega\subset\R^d$. We assume that $\U\subset\L$ is a continuous embedding. If $\cB$ has closed range the following are equivalent:
\begin{enumerate}
\item
$\cB:\U\to\L$ is an isomorphism, i.e., \eref{BN} holds for $\V=\L =\L'$, i.e.,
there exists a unique solution $u^*\in\U$, satisfying
\be
\label{1order}
b(u,v)=(\cB u)(v)= (f,v)_\Omega ,\quad v\in \L.
\ee
\item
For any $f\in\L$ there exists a unique minimizer
\be
\label{minu}
u^*=\argmin_{u\in\U}\|f- \cB u\|^2_\Omega,
\ee
that satisfies the following {\em error-residual relation}
\be
\label{err-res}
c_b\|u^*-u\|\le \|f- \cB u\|_\Omega \le C_b \|u^*-u\|,\quad u\in \U,
\ee
for some constants $0< c_b\le C_b<\infty$.
\item
The symmetric problem
\be
\label{normal}
B(u,v)= (\cB u,\cB v)_\Omega = (f,\cB v)_\Omega,\quad v\in \U,
\ee
is elliptic over $\U$ with coercivity and continuity constants $c_b^2, C_b^2$, respectively.
\end{enumerate}
\end{prop} 
\noindent
{\bf Proof:} Although the essence of the underlying reasoning is known we briefly sketch
the short arguments which are actually instructive.\\[2mm]
\underline{(i) $\Rightarrow$ (ii):} By the Babu\v{s}ka-Ne\v{c}as Theorem, $b(\cdot,\cdot)$
satisfies \eref{BN} for some positive constants $c_b, C_b$ and hence \eref{opbounds}.
Clearly the solution $u^*$ of \eref{1order} minimizes $\|f- \cB u\|_\Om$ and, due to its uniqueness, no other $\t u$ can minimize $\|f- \cB u\|_\L$. Moreover, since
$\|f-\cB u\|_\Omega = \|\cB(u^*-u)\|_\Omega$, \eref{err-res} follows from  \eref{opbounds}.\\[2mm]
\underline{(ii) $\Rightarrow$ (iii):} 
%$\frac 12 \|f-\cB u\|^2_\Omega$ is a continuous quadratic functional on $\U$ 
%whose values tend to infinity when $\|u\|_\U\to \infty$, there exists   a 
The unique minimizer $u^*$ of $\frac 12 \|f-\cB u\|^2_\Omega$ must
be a zero of the first variation of
$$
(f-\cB u,f-\cB u)_\Omega= \|f\|^2_\Omega + (\cB u,\cB u)_\Omega-2(f,\cB u)_\Omega,
$$
and hence a solution of the normal
equation \eref{normal}.
%the minimizer is  Uniqueness of the minimizer says that $u^*$ is the only solution to \eref{normal}. 
 Finally, every $u\in\U$ can be written as $u^*-u'$ for some $u'\in\U$. Thus,   \eref{err-res} says
$$
c_b\| u\|_\Omega \le \|\cB u\|,\quad \|\cB u\|_\Omega\le C_b\|u\|,
$$
which   implies (iii) by definition of $B(\cdot,\cdot)$.\\[2mm]
\underline{(iii) $\Rightarrow$ (i):} (iii) readily implies injectivity of $\cB$. Suppose, $\cB$ was not surjective. Since $\cB$ has closed range, there must exist a $g\in \L$ such that $(g,\cB u)=0$ for all $u\in\U$. On the other hand, Lax-Milgram ensures that there exists a $u_g\in \U\setminus \{0\}$ such that $(\cB u_g,\cB v)_\Omega= (g,\cB v)_\Omega$ for all $v\in \U$, hence
$c_b^2\|u_g\|_\U\le \|\cB u_g\|_\Omega = (g,u_g)_\Omega=0$, a contradiction. Hence, $\cB$ is surjective and thus bijective. By continuity of $B(\cdot,\cdot)$, $\cB$ is bounded. By the open mapping theorem, $\cB^{-1}:\L\to \U$ is also bounded, whence (i) follows from the Babu\v{s}ka-Ne\v{c}as Theorem.\hfill $\Box$

\begin{remark}
For convenience we have formulated \eref{opeq} or \eref{weak} as a {\em linear} operator equation, i.e.,
the trial space $\U$ already incorporates homogeneous essential boundary conditions. 
Depending on the specific problem at hand, inhomogeneous boundary conditions, can be either reduced
to such a case by modifying the right hand side or by enlarging the operator by trace components when these traces map into an $L_2$-space.
\end{remark}

Proposition \ref{prop:equiv} offers, in principle, three ways of numerically solving
$\cB u=f$. Since normal equations square condition numbers, the first order formulation
may be tempting. But since it is unsymmetric %(and not coercive) 
a discretization requires
for any trial space $\U_h\subset \U$ a judiciously chosen test space $\V_h\subset\L$. In the present
case the optimal test space (in the sense that it inherits the inf-sup constant from the
infinite dimensionall formulation) is simply $\V_h := \cB\U_h$ so that one ends up with
\eref{normal}.

%%%%%%%%%%%%%%%%%%%%%%%%%%%%%
\subsection{A Classical Example}\label{ssec:ex}
%%%%%%%%%%%%%%%%%%%%%%%%%%%%%%%
Let us briefly recall a by now familiar scenario where the above considerations apply,
namely the classical Poisson problem: for $g\in H^{-1}(\Omega)=(H^1_0(\Omega))'$ find
$u^*\in H^1_0(\Omega)$ such that
\be
\label{poisson}
( \nabla u^*,\nabla v)_\Omega = g(v),\quad v\in H^1_0(\Omega),
\ee
which is the weak formulation of $-\div(a\nabla u)=g$ in $\Omega$ with $u|_{\partial\Omega}=0$. 
Of course, this prohibits taking an $L_2$-type space as test space. But this changes 
when rewriting this second order equation as a {\em first-order system}. 
To that end, decomposing a general $g\in H^{-1}(\Omega)= (H^1_0(\Omega))'$ as $g= f_0 + \div \,\zeta$
for some $f_0\in L_2(\Omega)$, $\zeta\in L_2(\Omega;\R^d)$, the second order equation can be written as a system of first order equations
\be 
\label{system}
\sigma-  \nabla u = \zeta,\quad -\div\,\sigma= f_0.
\ee
Now one can seek $\bu^*=[\sigma^*,u^*]\in \U=H(\div;\Omega)\times H^1_0(\Omega)$ such
that for $\bbf=[\zeta,f_0]\in L_2(\Omega;\R^d)\times L_2(\Omega)=\L$
\be 
\begin{aligned}
\label{FOSLS}
(\cB \bu^*)(\bv)&= b(\bu^*,\bv) = b([\sigma^*,u^*],[\tau,v ])\\
&:= (\sigma^*-  \nabla u^*,\tau)_\Omega - (\div \,\sigma^*,v )_\Omega = (f_0,v )_\Omega + (\zeta,\tau)_\Omega,\\
&\qquad \forall\, \bv=[\tau,v ]\in \L,
\end{aligned}
\ee
i.e., 
\be
\label{B}
\cB \bu^* = {\sigma^*-  \nabla u^*\choose -\div \,\sigma^*} = {\zeta\choose f_0} =: \bbf.
\ee
In fact, it is shown in \cite{FOSLS} that (ii) in Proposition \ref{prop:equiv}
holds and hence, by Proposition \ref{prop:equiv},  isomorphy of $\cB: \U\to \L$,
as well as the minimum residual characterization.

The solutions of \eref{FOSLS} and the original problem \eref{poisson} agree in the following sense. For the solution $u^*\in H^1_0(\Omega)$ of \eref{poisson}   $\sigma^* := \nabla u^*,u^* =u^*$ solves  \eref{FOSLS}, and conversely, see
 \cite{FOSLS,Bochev}.
 
As a key point,  the setting considered in Section \ref{ssec:FOSLS} is viable when
dealing with first order systems. In fact, 
 very similar results have meanwhile been established for first order formulations of linear elasticity, for the Stokes system, as well as for
 space-time formulations of the heat or wave equation, \cite{FK,FGK} which motivates
 the discussion on the proposed more general level.

%%%%%%%%%%%%%%%%%%%%%%%%%%%%%%%%%
\subsection{Optimal Trial Norm}
%%%%%%%%%%%%%%%%%%%%%%%%%%%%%%%%%
We adhere to the above notation and    assume from now on that for a given Hilbert space $\U$, $\cB\in\cL(\U,\L)$ is
bijective so that all items in Proposition \ref{prop:equiv} are valid. We will make use of
 either of them whenever convenient. Thus, given $f\in \L$, there exists then
a unique $u^*\in\U$ such that
\be
\label{u*}
\cB u^* =f
\ee
in the sense of \eref{1order}.
We are interested in devising an {\em adaptive} procedure for approximately solving
\eref{u*}. Convergence analysis  will be greatly simplified when working with one specific (equivalent) norm on $\U$
which, of course, implies then convergence in any other equivalent norm. To be specific, due to 
\eref{opbounds},
\be
\label{optnorm}
\|u\|_\U := \|\cB u\|_\L,\quad u\in \U,
\ee
is an equivalent norm on $\U$, sometimes referred to as {\em optimal trial norm}, see e.g.
\cite{DGacta}. Specifically, well-posedness of \eref{FOSLS} is equivalent to saying that in this case
$$
\begin{aligned}
\|[\sigma,u]\|^2_{H(\div;\Omega)\times H^1_0(\Om)}&= \|\sigma\|_{L_2(\Omega;\R^d)}^2+\|\div\,\sigma\|^2_{L_2(\Om)}+ \|\nabla u\|^2_{L_2(\Om)} \\
& \eqsim \|\cB[\sigma,u]\|_{L_2(\Om;\R^{d+1})} \\
&  = \|\sigma-a\nabla u\|^2_{L_2(\Om;\R^d)}+ \|\div\,\sigma\|^2_{L_2(\Om)}, 
  \end{aligned}
$$
holds for all $ [\sigma,u]\in H(\div;\Omega)\times H^1_0(\Om).
$
 The reason
for calling \eref{optnorm} optimal  is that for this norm one has by definition $\|\cB\|_{\cL(\U,\L)}=1$ while also
$$
\inf_{u\in\U}\sup_{v\in\L}\frac{b(u,v)}{\|u\|_\U\|v\|_\L}=\inf_{u\in\U}\sup_{v\in\L}\frac{(\cB u,v)_\Omega}{\|\cB u\|_\L \|v\|_\L}=\inf_{u\in\U}\frac{\|\cB u\|_\L}{\|\cB u\|_\L}=1,
$$
i.e., $c_b=C_b=1$. Thus, $\cB$ is even an isometry (condition number equals  one) and
we have a tight error-residual relation
\be
\label{tight}
\|u^*-u\|_\U = \|f- \cB u\|_\L,\quad u\in\U.
\ee

While the equivalence of the ``optimal trial norm'' \eqref{optnorm} with an ``original generic trial''
$\|\cdot\|$ on $\U$ - in the example from Section \ref{ssec:ex} $\|[\sigma,u]\|^2= \|\sigma\|^2_{H({\rm div};\Omega)} + \|u\|^2_{H^1(\Omega)}$ - follows from well-posedness in the above sense, corresponding
equivalence constants may strongly depend on the range of coefficients in $\cB$ or on near-degeneracies of 
the domain  $\Omega$ and may suggest specific scaling weights in different terms of the $\L$-residual. However, that does not affect the principal question of error reduction of cycles like \eqref{cycle}. Again, on the current level of generality, we cannot and don't have to address this issue in what follows.
%%%%%%%%%%%%%%%%%%%%%%%%
\newcommand{\mT}{\mathfrak{T}}
%%%%%%%%%%%%%%%%%%%%%%%%%
\section{An Idealized Setting}\label{sec:prel}

\subsection{Core Ingredients of an Adaptive Cycle}\label{ssec:ingredients}
%%%%%%%%%%%%%%%%%%%%%%%%%%%%%%%%%%%%%%

We denote by $\cT$ a generic partition of a given computational domain $\Omega$, subject
to certain admissibility constraints, common in finite element methods. The precise specification of such conditions may depend on the underlying PDE model.
They may require geometric conformity, meaning the intersection of the closure of adjacent cells must be a common facet of both cells. Admissibility may also allow for hanging nodes whose number per element facet may, however, be uniformly bounded.
When viewing $\cT$ as an element of a collection $\mT $ of partitions,
a typical further constraint is that the $\cT$ in $\mT$  are uniformly shape regular.
We focus   on {\em nested hierarchies} $\mT=\{\cT_j\}_{j\in\N_0}$ 
%whose elements are {\em nested}
\be
\label{nested}
\cT_{ 0}\prec \cT_{ 1}\prec \cdots\prec \cT_{n}\prec \cdots
\ee
where $\cT_j\prec\cT_{j+1}$ means that $\cT_{j+1}$ arises from $\cT_j$ by subdividing
some cells in $\cT_j$ according to a specific rule, possibly followed by closing refinements so as     to preserve {\em admissibility},
see below for further comments in this regard. We use $\cT\in \mT$ synominously with saying $\cT$
is $\mT$-admissible. We always assume that the specific partition $\cT_0^j$,   obtained by a $j$-fold subdivision 
of {\em all cells} in $\cT_0$, is $\mT$-admissible. A prototype example is {\em newest verstex bisection (or longest-edge-bisection)}, see \cite{BDD,S}.

With each $\cT\in \mT $ we associate a finite element trial space $\U_\cT \subset\U$.
In the example of Section \ref{ssec:ex} one would use $H(\div)$- and $H^1$-conforming finite elements,
see e.g. \cite{Braess,BBF}.
This gives then rise to a hierarchy of nested finite elements spaces 
$$
\U_{\cT_0}\subset \U_{\cT_1}\subset \cdots \subset \U_{\cT_j}\subset \cdots.
$$
 For each $\cT\in\mT$ 
 let
\be
\label{uT}
u_\cT=\argmin_{u\in \U_\cT}\frac 12\|f- \cB u\|^2_\Omega,
\ee
which (by Proposition \ref{prop:equiv}) is the solution of
\be
\label{normalT}
(\cB u_\cT,\cB v)_\Omega= (f,\cB v)_\Omega,\quad v\in \U_\cT.
\ee
For $\cT=\cT_j$ we briefly write at times $u_j:= u_{\cT_j}$.

Denoting for any closed subspace $\wh\U\subset\U$ by $P_{\wh\U}$ the $\U$-orthogonal
projection with respect to the optimal trial norm $\|\cdot\|_\U= \|\cB \cdot\|_\Om$, 
\eref{normalT} means in particular
\be
\label{OU}
u_\cT = P_{\U_\cT} u^*,
\ee
where $u^*$ consistently denotes the exact solution to \eref{u*}.
We formalize this by the routine \\[2mm]
$\textsc{Solve}[f,\cT]\to u_\cT= u_\cT(f)$\\[2mm]
that generates the {\em exact} solution of the discrete system \eref{normalT}.

This is a first instance of a strong ``idealization'' that will be adopted in a first version of
an adaptive algorithm. In practice, neither can we solve \eqref{normalT} exactly  nor can
we generally compute the matrix and load vector in \eqref{normalT} exactly. It would be more realistic to
talk about $\textsc{Solve}[f,\cT,\delta]$ meaning that we only insist on solving the system within a prescribed tolerance $\delta>0$, see e.g. \cite{robfocm}. We prefer to suppress such aspects for the moment for the sake of simplicity and clarity of exposition. Instead we   refer to corresponding
existing strategies e.g. to \cite{robfocm}, that apply in the present context as well.
%We will however take up this aspect later with focus on perturbing the load $f$ because 
%it turns out to directly related  to the essence of convergence issues.  

We proceed formally specifying the ingredients of \eref{cycle} as ``routines'', beginning with marking the cells to be
refined.\\ 

\noindent
$\textsc{Mark}[\cT,u_\cT,f,\theta]\to \cK_\theta(\cT,u_\cT,f)= \cK_\theta$\\[2mm]
 We adopt in what follows {\em D\"orfler marking} as the way to determine a subset $\cK\subset \cT$ that should be refined. Namely
fix $\theta\in (0,1)$ and   let $\cK_\theta=\cK_\theta(\cT)$ be a minimal subset of $\cT$
such that
\be
\label{Ktheta}
\sum_{K\in \cK_\theta} \|f-\cB u_\cT\|^2_K\ge \theta\|f-\cB u_\cT\|^2_{\Omega}.
\ee
Again, at this point this is an idealization because the element-residuals $\|f-\cB u_\cT\|^2_K$
are computed exactly. 
We call $\cK_\theta$ to be {\em minimal} if, up to a fixed multiplicative factor,
$\#\cK_\theta(\cT)$ is minimal among all collections satisfying \eref{Ktheta}.

Clearly, abbreviating 
$$
\Om_\cK:= {\rm interior}\,\Big(\bigcup_{K\in \cK}\bar K\Big),
$$
\eref{Ktheta}
implies 
\be
\label{bulk}
 \|f-\cB u_\cT\|^2_{\Om_\cK }\ge \theta\|f-\cB u_\cT\|^2_{\Omega}.
 \ee

Once   $\cK_\theta$ has been identified, the subsequent {\em mesh refinement}
is
 based on several further ingredients.   
\\[2mm]
$\textsc{Subd}[\cT,\cK]\to \cT_{\cK}$\\[2mm]%
 subdivides for a fixed number $\ell\ge 1$ times all cells in $\cK \subset\cT$
 according to the rule underlying   $\mT$. Since $\ell$ remains fixed we suppress any notational reference.

The  partition $\cT\lor \cT_{ \cK}$ may not be admissible yet and may need to be {\em closed} in $\mT$
to restore admissibility.  %  
  
The routine \\

\noindent
$\textsc{Refine}[\cT,\cK]\to \cTr(\cT,\cK)$,\\[2mm]
generates the smallest admissible refinement of $\cT$ that contains $\textsc{Subd}[\cT,\cK]$.
In particular, this gives rise to the {\em refinement region}    
 $$
 \Or=\Or(\cT,\cK) = {\rm interior}\,\Big(\bigcup\big\{ K:K\in \cTr\setminus\cT  \big\}\Big).
 $$
which is the union of subdivided cells in $\cT$. Whenever $\cT,\cK$ are clear from the context we briefly write $\Or $.

%%%%%%%%%%%%%%%%%%%%%%%%%%%%%%%%%%%%%
\subsection{  \textsc{ALSM}  - an Idealized Version}\label{ssec:ideal}
%%%%%%%%%%%%%%%%%%%%%%%%%%%%%%%%%%%%%%%%%%%%%%%%%%%%%%%%%%%
We adhere to the above idealizing assumptions and abbreviate at times for convenience   for any $v\in \U_\cT$
$$
\cR_\cT( v,f):= f- \cB v, \,\,\mbox{ and in particular }\,\, { \cR_\cT := \cR_\cT(u_\cT,f)}.
$$
To be concise we formalize  \eqref{cycle} as follows:\\[2mm]
\noindent
$\textsc{Alsm-Ideal}[\cB,f ,\theta ,\e] \to \cT_\e,\, u_\e\in \U_{\cT_\e}$: 
\begin{itemize}
\item[(A0)]
Choose $\cT_0$, fix $\theta, \in (0,1)$,   $\e>0$;\\  
set $0\to \bar u$, $\cT_0\to \cT$;  
\item[(A1)]
compute $ \|f- \cB \bar u\|_\Om$;\\
while $\|f- \cB \bar u\|_\Om>\e$ apply $\textsc{Mark}[\cT,f,\bar u ,\theta]\to \cK_\theta$;\\
else, output $\cT_\e= \cT$, $u_\e=\bar u$;
\item[(A2)] 
apply $\textsc{Refine}[\cT,\cK_\theta  ]\to  \cTr  $;
\item[(A3)]
$\textsc{Solve}[f, \cTr]\to   u_\cTr$;\\ 
set $\cTr\to \cT, u_\cTr\to \bar u$ and
go to (A1).
 \end{itemize}
\medskip

Obiously, if the scheme terminates the output satisfies $\|u^*-u_\e\|_\U= \|f- \cB u_{\cT_\e}\|_\Om\le \e$.

Understanding whether or under which circumstances  \textsc{Alsm-Ideal} terminates and
how   a corresponding convergence can be quantified, will pave the way to formulating and analyzing subsequent
more practical versions.

%%%%%%%%%%%%%%%%%%%%%%%%%%%%%%%%%%
\subsection{What Can be Expected?}
%%%%%%%%%%%%%%%%%%%%%%%%%%%%%%%%%%%%%%%
\subsubsection{Error reduction and the ``adaptive update''}
%%%%%%%%%%%%%%%%%%%%%%%%%%
 We begin with the (very well-known) Galerkin orthogonality
 and its simple manifestation in the current setting
 \be
\label{Udec}
\|u^*- u_\cTr\|^2_\U = \|u^*-u_\cT\|^2_\U - \|u_\cTr-u_\cT\|^2_\U\ge 0.
\ee
Hence for some $s\in [0,1]$ one must have  
\be
\label{smust}
\|\wT\|^2_\U=
s\|u^*- u_\cT\|^2_\U \quad \Leftrightarrow  \quad\|u^*-u_\cTr\|^2_\U= (1-s)\|u^*- u_\cT\|^2_\U.
\ee
We record the following immediate consequence of \eref{smust}.
 
\begin{rem}
\label{rem:iff}
Adhering to the above notation, one has
\be
\label{equiv}
\begin{array}{c}
\|u^*- u_\cTr\|_\U \le \rho\|u^*- u_\cT\|_\U\\
 \Leftrightarrow\\
  \|u_\cT-u_\cTr\|_\U \ge
\sqrt{1-\rho^2}\|u^*-u_\cT\|_\U.
\end{array}
\ee
In other words, the better the ``adaptive update'' 
\be
\label{udec}
\wT= \wT(f,\cTr):=  u_\cTr- u_\cT
\ee
aligns with $u^*-u_\cT$, the better the error reduction.
\end{rem}

We'll make repeated use of the relations 
\be
\label{worth}
\wT \perp_\U u^*-u_\cTr,\quad \wT \perp_\U \U_\cT,
\ee
which follow from the fact that, as solutions to \eref{normalT},  $u_\cT$, $u_\cTr$, are the $\U$-orthogonal projections of $u^*$ to $\U_\cT, u_\cTr$, respectively, see \eref{OU}.

The key is to understand under which circumstances $\|\wT\|_\U=\|\cB \wT\|_\Om$ captures a ``substantial part'' of the residual $\|f- \cB u_\cT\|_\Om = \|u^*- u_\cT\|_\U$.
 To that end,  the second relation in \eref{worth} suggests    decomposing the refined space $\U_\cTr=\U_{\cT,\cTr(\cK)}$ as
 $$
\U_\cTr= \U_\cT \oplus_\U \W
 $$ 
where
\be
\label{Wdef}
\W=\W(\cT;\cTr ):= P_{\U^\perp_\cT} \U_{\cTr}:= (I - P_{\U_\cT})\U_\cTr\perp_\U \U_\cT. 
\ee
Notice further that, since $\U_\cT\subseteq \U_\cTr$, so that $P_{\U_\cT}$ and $P_{\U_\cTr}$ commute, i.e.,
$$
(I - P_{\U_\cT})\U_\cTr= (I - P_{\U_\cT})P_{\U_\cTr}\U = P_{\U_\cTr}(I - P_{\U_\cT})P_{\U_\cTr}\U
=(P_{\U_\cTr}- P_{\U_\cT})\U,
$$
we see that $P_{\U_\cTr}- P_{\U_\cT}$ is also a projector. Hence,   we have 
\be
\label{Wl2}
\W = (P_{\U_\cTr}-P_{\U_\cT})\U=: P_\W \U.
\ee
%%%%%%%%%%%%%%%%%%%%%
\subsubsection{Two Ways   to look at the adaptive update $\wT$}
%%%%%%%%%%%%%%%%%%%%%%%%

First, it follows that 
\be
\label{wTrep}
\wT = (P_{\U_\cTr}-P_{\U_\cT})u^*=  (P_{\U_\cTr}-P_{\U_\cT})(I-P_{\U_\cT})u^*= P_\W u^*,
\ee
that is, $\wT$ is the $\U$-projection of $u^*$ to $\W$.

On the other hand, it is at times convenient to  interpret the update $\wT$ as the solution to a variational problem, posed in $\W$.
 \begin{rem}
\label{rem:w0}
The   correction $\wT\in \W\subset\U_\cTr$ solves
\be
\label{w}
(\cB \wT,\cB v)_\Omega = (\cR_\cT  ,\cB v)_\Omega,\quad v\in \U_\cTr,
\ee
so that
(see \eref{Wdef} and the second relation in \eref{worth})
\be
\label{w2}
\begin{aligned}
&(\cB \wT,\cB v_\cT)_\Omega =0, & \quad v_\cT\in\U_\cT, %\\
\end{aligned}
\ee
and therefore
\be
\label{Wproj}
(\cB\wT,\cB w)_\Om = (f,\cB w)_\Om,\quad w\in \W.
\ee
In particular
\be
\label{wportion}
\|\wT\|_\U^2= (\cR_\cT,\cB \wT)_\Om = (f,\cB \wT)_\Om = (u^*,\wT)_\U. %= (u^*-u_\cT,\wT)_\U.
\ee 
\end{rem} 

Thus, how well $\cB\wT$ aligns with the residual in $\L$ is equivalent to its alignment with
just the data $f$.
\begin{rem}
\label{rem:many}
It follows from \eref{wportion} that whenever $f$ happens to be $\L$-orthogonal to $\cB\wT$, there is no
error reduction at all. Since $\W$ is finite dimensional, this already shows that without any further condition on $f\in\L$, one cannot exclude that   for
any finite subdivision depth $\ell\in\N$ in \textsc{Subdiv}, the ``extreme''
case
\be
\label{extreme}
  f  \perp_\L \cB \W
\ee
could happen preventing any error reduction. By the same token, error reduction can be arbitrarily weak
when the angle between $f$ and $\cB\W$ is small which may be the case even when $f$ is pointwise very smooth. 
\end{rem} 
 
\subsubsection{Data-Visibility}\label{ssec:osc}
%%%%%%%%%%%%%%%%%%%%%%%%%%%%%%%%%%%%%%%%%%%%%%
\newcommand{\Uinfty}{{\U_{\cT}(\Or)}}
%%%%%%%%%%%%%%%%%%%%%%%%%%%%%%%
Intuitively, one expects that \eref{extreme} won't happen, or that the angle between $f$ and $\W$
eventually gets smaller, if $\ell$ - hence $\W$ - is sufficiently large. To put this thought
into proper perspective, we introduce certain subspaces of $\U$, associated with
subdomains $\Om'\subset \Om$ (typically unions of cells in some partition from $\mT$).
Specifically,
 we denote   by $\U(\Om')$ the following closed subspace of 
$\U=\U(\Om)$:
\be 
\label{Uloc}
\U(\Om'):= \big\{ u\in \U: \chi_{\Om'}u\in \U\big\}.
\ee
In other words, $\U(\Om')$ is comprised of those $u\in\U$ whose trivial extension by zero to $\Om\setminus
\Om'$ still belongs to $\U$. We call
a subdomain $\Om'\subset \Om$   {\em regular} if \eqref{normal} is well-posed on $\U(\Om')$
 for any $g\in \L(\Om')$, e.g. for $g=f|_{\Om'}$ when $f\in \L$. In particular, simply connected 
 unions of cells in $\cT$ or disjoint unions of such connected subdomains are regular.

In addition, for any partitions $\cT\prec\cT'$ in $\mT$ we will frequently work with the subspaces
\be
\label{subspaces}
\begin{array}{c}
\U(\Om';\cT'):= \U(\Om')\cap \U_{\cT'}, \quad \U_\cT(\Om'):= \U_\cT +\U(\Om'),\\ 
\U_\cT(\Om';\cT')= \U_\cT + \U(\Om';\cT').
\end{array}
\ee
Specifically, we can write in these terms
\be
\label{write}
\U_\cTr = \U_\cT\oplus_\U \W = \U_\cT + \U(\Or;\cTr)= \U_\cT(\Or;\cTr).
\ee 
which is obviously contained in
\be
\label{Uinfty}
 \U_\cT + \U(\Or)= \U_\cT(\Or),
\ee
This latter space can be formally viewed as resulting form letting the refinement depth $\ell$ tend to infinity on $\Or$.
Now, in analogy to \eref{Wdef}, let
\be
\label{Wi}
\Wi := (I- P_\cT)\Uinfty = (P_{\Uinfty }-P_{\U_\cT})\U= (I- P_\cT)\U(\Or),
\ee
which obviously contains $\W$ for every $\ell\in \N$.

We are now prepared, to introduce a notion that quantifies how 
``visible'' $f$
- hence the residual $f-\cB u_\cT$ - 
is by $\cB \W$:  
\be
\label{inv}
{ \osc_\cT(f,\W )}:= \sup_{u\in \U}\frac{(f,\cB (P_\Wi- P_{\W})u)_\Omega}{\|\cB u\|_\Omega}={ \|(P_\Wi- P_{\W})'\cB'f\|_{{\U}'}}.
\ee
Obviously, $\osc_\cT(f,\W )$ is large when $(P_\W'\cB'f)(w)=(f,\cB w)_\Om$ is small for each $w\in\W$.

There is an obvious alternate characerization of $ \osc_\cT(f,\W )$ as a projection error.
\begin{lemma}
\label{lem:rep}
One has
\be
\label{oscrep}
\osc_\cT(f,\W) = \|(P_{\Wi}- P_{\W})u^*\|_\U = \frac{(\cR_\cT,\cB (P_\Wi- P_{\W})u^*)_\Omega}{\|(P_{\Wi}- P_{\W})u^*\|_\U}.
\ee
\end{lemma}
\noindent
{\bf Proof:} For each $u\in \U$ one has
\begin{align*}
\frac{(f,\cB (P_\Wi- P_{\W})u)_\Omega}{\|u\|_\U} &= \frac{(u^*,(P_\Wi- P_{\W})u)_\U}{\|u\|_\U}\\
&\le  \frac{(u^*,(P_\Wi- P_{\W})u)_\U}{\|(P_\Wi- P_{\W})u\|_\U} ,
\end{align*}
which shows that $\osc_\cT(f,\W)\le \|(P_\Wi- P_{\W})u^*\|_U$. On the other hand,
by similar arguments
\begin{align*}
\sup_{u\in\U}\frac{(f,\cB (P_\Wi- P_{\W})u)_\Omega}{\|u\|_\U}&\ge \frac{(f,\cB (P_\Wi- P_{\W})u^*)_\Omega}{\|(P_\Wi- P_{\W})u^*\|_\U}\\
&= \|(P_\Wi- P_{\W})u^*\|_U,
\end{align*}
proving the claim since $(f,\cB (P_\Wi- P_{\W})u)_\Omega=(\cR_\cT,\cB (P_\Wi- P_{\W})u)_\Omega$.
\hfill $\Box$\\

The quantity $\osc_\cT(f,\W )$ formally resembles  what one may call {\em data-oscillation} of $f$.
It is, however, closely tied to the current refinement region $\Or$.
 
%%%%%%%%%%%%%%%%%%%%%%%%%%%%%%%%%%%%%%%
\subsection{An Idealized Convergence Result}\label{ssec;idconv}
%%%%%%%%%%%%%%%%%%%%%%%%%%%%%%%%%%%%%%%%%%%%%%
The scheme \textsc{Alsm-Ideal} from Section \ref{ssec:ideal} gives rise to a sequence of quantities which are denoted throughout the remainder 
of this section as follows:
\be
\label{notation}
\mbox{\begin{tabular}{l}
-  $(\cT_j)_{j\in\N_0}$: hierarchy of nested partitions;\\

-  $u_j=
 u_{\cT_j}\in \U_j:= \U_{\cT_j}$:  associated discrete solutions,\\ 
 
 \hspace*{3mm}($\cT_j=\cT\to \cTr =\cT_{j+1}$), $j\in \N_0$;
\\
 
 - $\W_j:= \W(\cT_j,\cT_{j+1})$, see \eqref{Wdef}.
 \end{tabular}
 }
 \ee
\begin{theorem}
\label{thm:main1}
We adhere to previous notation.\\[2mm]
(a) Given $\theta\in(0,1)$, 
then,   the output $u_\cTr$ of a single cycle of \textsc{Alsm-Ideal}
satisfies
\be
\label{reliable}
\|u^*- u_\cTr\|_\U^2 \le (1-\theta)^{1/2} \|u^*-u_\cT\|^2_\U +\osc_\cT(f,\W)^2,
\ee
\small
(b) If there exists a constant $C_0<\infty$ for which
\be
\label{ell2}
\sum_{j=k}^\infty \osc_{\cT_j}(\W_j;f)^2 \le C_0\|\cR_{\cT_k}\|^2_\Om,\quad k\in \N_0,
\ee
holds uniformly in $k\in\N_0$, then there exists a fixed $m\in \N$ and a $\rho\in (0,1)$,
depending only on $\theta$, $\ell$, and $C_0$
such that
\be
\label{linear1}
\|u^*- u_{j+m}\|_\U\le \rho \|u^*-u_j\|_\U,\quad j\in \N_0.
\ee
\end{theorem}
%

%Note that, while $\rho, m$ may depend on $\ell,\theta$, a linear error reduction \eqref{linear1} occurs regardless of their specific
%choice.  
Of course, \eqref{linear1}
can be rephrased as: there exists a constant $C_1$, depending only on $m,\rho$  such that for $\t\rho=\rho^{1/m}<1$
\be
\label{linear2}
\|u^*- u_n\|_\U\le C_1\t\rho^n\|f\|_\Om,\quad n\in\N.
\ee

\begin{rem}
\label{rem:usual}
The  linear error reduction, up to the data-visibility term, \eref{reliable} is in this form
unsatisfactory because it is not clear yet that $\osc_\cT(f,\W)$ is small when $f$ is smooth. 
So the question is under which assumptions on $f$ is $\osc_\cT(f,\W)$ dominated by the residual,
which motivates  the subsequent reasoning.
\end{rem}

%%%%%%%%%%%%%%%%%%%%%%%%%
\subsection{Proof of Theorem \ref{thm:main1}}\label{ssec:Proof1}
%%%%%%%%%%%%%%%%%%%%%%%%%%%%%%%%%%%%
We proceed collecting a few technical prerequisits.
%%%%%%%%%%%%%%%%%%%%%%%%%%%%%%
\subsubsection{Further elementary properties of $\osc_\cT(f,\W)$}
%%%%%%%%%%%%%%%%%%%%%%%%%%%%%%%%%%%%%%%%%%%%%
Note that for any $w\in \Wi$ one has $(\cR_\cT,\cB w)_\Om = (f,\cB w)_\Om= (u^*,w)_\U$. 
The definition  \eref{Wi} of $\Wi$ suggests then considering the ``infinite-counterpart'' of $\wT$, namely the
projection of $u^*$ to $\Wi$
\be
\label{wT*}
\wT^*:=(I-P_{\U_\cT})P_\Uinfty u^*=P_\Uinfty (I-P_{\U_\cT})u^*=P_\Wi u^* \in \Wi.
\ee
This gives rise to the following analog to Remark \ref{rem:w0}
\begin{rem}
\label{rem:wT*}
Adhering to the above notation,  
$\wT^*$ solves the  variational problem in $\Uinfty$: 
\be
\label{w*}
 (\cB  \wT^*,\cB v)_\Om  = (\cR_\cT,\cB v)_\Om
 ,\quad v\in \Uinfty.
\ee
Moreover,
one has
\be
\label{wT*2}
\|\wT^*\|_\U^2=\|\cB \wT^*\|^2_\Om = (\cR_\cT,\cB \wT^*)_\Om =\|P_{\Wi}u^*\|_\U^2,
\ee
i.e. the $\U$-energy of the ``infinitely refined  update agrees with the one of the projection
of $u^*$ to $\Wi$.
\end{rem}
\noindent
{\bf Proof:} Regarding \eref{w*},
$$
\begin{aligned}
(\cR_\cT,\cB v)_\Om &= (u^*-u_\cT,v)_\U = (u^*,(I- P_{\U_\cT})P_\Uinfty v)_\U \stackrel{\eref{Wi}}{= }(u^*,P_\Wi v)_\U\\
&= (P_\Wi u^*,v)_\U = (\cB \wT^*,\cB v)_\Om,\quad v\in \Uinfty,
\end{aligned}
$$
which also gives \eref{wT*2} upon taking $v=\wT^*\in \Wi\subset \U_\cT(\Or)$.
\hfill $\Box$\\

\begin{rem}
\label{rem:f'}
In particular, one has 
\be
\label{large}
f\perp_\L \cB \W \quad\Leftrightarrow\quad \osc_\cT(f,\W) =\|P_\Wi u^*\|_{\U}.
\ee 
indicating, in view of Remark \ref{rem:many}, how $\osc_\cT(f,\W)$, relates to 
error reduction.
\end{rem}

%%%%%%%%%%%%%%%%%%%%%%%%
\subsubsection{Data visibility and residuals}
%%%%%%%%%%%%%%%%%%%%%%%%%%%%%%%%%%
The quantity $\osc_\cT(f,\W)$ always bounds the current residual, hence the error, from below.

\begin{rem}
 \label{rem:below}
 One has
 \be
 \label{refbelow}
 \osc_\cT(f,\W) \le \|\cR_\cT\|_\Om.  
 \ee
  \end{rem}
 \noindent
 To see this, recall from Lemma \ref{lem:rep} that
 \begin{align*}
 \osc_\cT(f,\W)= \|(P_{\Wi}-P_\W)u^*\|_\U &= \frac{(\cR_\cT,\cB (I-P_\W)P_{\Wi}u^*)_\Om}{
 \|(P_{\Wi}-P_\W)u^*\|_\U}\\
 &\le \|\cR_\cT\|_\Om ,
  \end{align*}
 whence the assertion follows.\hfill $\Box$\\

 The key to strict error reduction  lies  in tightening \eref{refbelow}.
 
\begin{lemma}
\label{lem:supp}
Every $z\in \Uinfty^\perp$ must satisfy
\be
\label{suppz}
{\rm supp}\,\cB z \subseteq \overline{\Or^c}.
\ee
\end{lemma}
\noindent
{\bf Proof:} Every $\t u\in \Uinfty$ can be written as $\t u= v_\cT + u_{\Or}$, where
$v_\cT\in \U_\cT$ and $u_{\Or}\in \U(\Or)$. Since $z\in \Uinfty^\perp$,  meaning that
 $z\perp \Wi$ and $z\perp \U_\cT$, one obtains
$$
0=(z,\t u)_\U = (z,u_{\Or})_\U= (\cB z,\cB u_{\Or})_\Om.
$$
Since $\Or$ is a regular polygonal domain $\cB$ maps $\U(\Or)$ onto $\L(\Or)$,  which confirms the claim. 
\hfill$\Box$\\

We are now prepared to establish the following key facts.
\begin{prop}
\label{lem:key}
Given $\theta\in(0,1)$, let $\xi=\xi(\theta)$ 
  be defined by  
\be
  \label{xi}
 \xi:= \big(1-(1-\theta)^{1/2}\big)^{1/2}.
 \ee
 Then
 \be
\label{lowerw*}
\|\cB \wT^*\|_\Om \ge \xi \|\cR_\cT\|_\Om ,
\ee
while
 \be
 \label{osc4}
 \|\cB \wT\|^2_\Om\ge  \xi^2 \|\cR_\cT\|^2_\Om- \osc_\cT(f,\W)^2.
 \ee
\end{prop}
Note that, while the alignment of the ``infinite increment'' $\wT^*$ with $u^*-u_\cT$ depends only
 on the bulk constant $\theta$, the alignment of the adaptive update $\wT$ with $u^*-u_\cT$,   which is our primary concern,   depends also on $\osc_\cT(f,\W)$.\\
 
\noindent
{\bf Proof:} Abbreviate  $e_\cT^*:= u^*- u_\cT$ and
recall from \eref{wT*} that 
\be
\label{wTreps}
\wT^*:= (I-P_{\U_\cT})P_{\Uinfty }u^*=
 P_{\Uinfty }  e^*_\cT \in \Wi.
\ee
Thus, 
\be
\label{proj2}
e^*_\cT-\wT^*= (I- P_{\Uinfty})e_\cT^*\in \Uinfty^\perp,
\ee
which gives
\be
\label{z}
\|e_\cT^*-\wT^*\|_\U =\sup_{z\in \Uinfty^\perp}\frac{(e_\cT^*-\wT^*,z)_\U}{\|z\|_\U}
= \sup_{z\in \Uinfty^\perp}\frac{(e_\cT^*,z)_\U}{\|z\|_\U}.
\ee
Invoking Lemma \ref{lem:supp}, yields
\begin{align}
\label{g*bound}
\|e_\cT^*-\wT^*\|_\U&= \sup_{z\in \UTr^\perp}\frac{(\cB e_\cT^*,\cB z)_{\Or^c}}{\|z\|_\U}
\le \|\cB e_\cT^*\|_{\Or^c}= \|\cR_\cT\|_{\Or^c}\nonumber\\
&\le (1-\theta)^{1/2}\|\cR_\cT\|_\Om= (1-\theta)^{1/2}\|e_\cT^*\|_\U.
\end{align}
Moreover, since $\wT^*\in \Uinfty$, and thus
$$
(e_\cT^*-\wT^*,\wT^*)_\U= ((I-P_{\Uinfty})e_\cT^*,\wT^*)_\U= (e^*_\cT,(I-P_{\Uinfty})\wT^*)_\U=0,
$$
we conclude, upon using that $\cB e^*_\cT = \cR_\cT$,
\be
\label{ew}
\begin{aligned}
(e_\cT^*,\wT^*)_\U&=(\cB e_\cT^*,\cB \wT^*)_\Om =(\cR_\cT,\cB \wT^*)_\Om 
 \stackrel{\eref{w*}}{= }\|\cB \wT^*\|^2_\Om.
\end{aligned}
\ee
On the other hand,
\begin{align*}
(e_\cT^*,\wT^*)_\U&= \|e^*_\cT\|_\U+ (e^*_\cT,\wT^*- e^*_\cT)_\U= \|\cR_\cT\|^2_\Om +(\cR_\cT,\cB(\wT^*- e^*_\cT))_\Om\\
&\ge\|\cR_\cT\|^2_\Om - \|\cR_\cT\|_\Om \|e_\cT^*- \wT^*\|_\U\nonumber\\
&\ge (1- (1-\theta)^{1/2})\|\cR_\cT\|^2_\Om,
\end{align*}
where we have used \eref{g*bound} in the last step.
This latter inequality together with \eref{ew}   provides
\begin{align*}
\|\cB \wT^*\|^2_\Om &\ge 
  (1- (1-\theta)^{1/2})\|\cR_\cT\|^2_\Om ,
\end{align*}
which, in view of the definition, 
 $\xi:= (1- (1-\theta)^{1/2})^{1/2}$ from \eref{xi}, confirms \eref{lowerw*}.

Regarding \eqref{osc4},   observe that
\begin{align}
 \label{osc40}
\|\cB \wT\|^2_\Om& \stackrel{\eref{wportion}}{=}(\cR_\cT,\cB\wT)_\Omega 
 = (\cR_\cT,\cB \wT^*)_\Omega -  (f, \cB(\wT^*- \wT))_\Om.
 \end{align}
Now, recall from \eqref{oscrep} in Lemma \ref{lem:rep} that
$
(f, \cB(\wT^*- \wT))_\Om = 
     \osc_\cT(f,\W)^2$.
Combining this with \eref{osc40}, yields
\begin{align*}
\|\cB \wT\|^2_\Om& \stackrel{\eref{wportion}}{=}(\cR_\cT,\wT)_\Omega 
 = (\cR_\cT,\cB \wT^*)_\Omega - \osc_\cT(f,\W)^2  \nonumber\\
&\stackrel{\eref{wT*2}}{=} \|\cB \wT^*\|^2_\Om -  \osc_\cT(f,\W)^2\nonumber\\
&\stackrel{\eref{lowerw*}}{\ge} \xi^2 \|\cR_\cT\|^2_\Om- \osc_\cT(f,\W)^2.
\end{align*}
as claimed. \hfill $\Box$
%%%%%%%%%%%%%%%%%%%%%%%%
\subsubsection{Proof of Theorem \ref{thm:main1}, (a)}
%%%%%%%%%%%%%%%%%%%%%%%%%%%%%%%%%%%%%%%%%%%%%%%
Note that, in view of \eqref{osc4},  
\begin{align*}
\|u^*-u_\cTr\|^2_\U &= \|u^*- u_\cT\|^2_\U - \|\wT\|^2_\U\\
&\le \|u^*- u_\cT\|^2_\U -\xi^2\|u^*-u_\cT\|^2_\U + \osc_\cT(f,\W),
\end{align*}
 part (a)
of Theorem \ref{thm:main1} follows from \eref{xi}.\\

%%%%%%%%%%%%%%%%%%%%%%%%
\subsubsection{Proof of Theorem \ref{thm:main1}, (b)}
%%%%%%%%%%%%%%%%%%%%%%%%%%%%%%%%%%%%%%%%%%%%%%%
To complete the proof of Theorem \ref{thm:main1} we need one further technical preliminary.
 
\begin{lemma}
\label{lem:xiC0}
Given $\xi$, $C_0$ according to \eref{xi}, \eref{ell2}, respectively.
Then there exist a $\rho<1$ and some $m\in\N$ such that
\be
\label{solution}
\rho^2\Big\{\frac{\xi^2(1- \rho^{2(m-1)})}{1-\rho^2} +1\Big\}\ge C_0+1.
\ee
\end{lemma}
\noindent
{\bf Proof:}
Define
$$
G(\rho,\xi):= \rho^2\Big\{\frac{\xi^2 }{2(1-\rho^2)} +1\Big\}.
$$
Clearly, $G(1,\xi)=\infty$ and $G(0,\xi)=0$. By monotonicity, there exists a unique
$\rho=\rho(\xi)$ such that $G(\rho,\xi)= C_0+1$. For this $\rho$ pick $m\in\N$ such that
$\rho^{2(m-1)}<1/2$. Then
$$
\rho^2\Big\{\frac{\xi^2(1- \rho^{2(m-1)})}{1-\rho^2} +1\Big\}\ge G(\rho,\xi)= C_0+1,
$$
as claimed.\hfill $\Box$\\

\noindent
{\bf Proof of Theorem \ref{thm:main1}, (b):} 
We adhere to the notation \eref{notation}. That is,  let $\cT=\cT_k$, $\cT_{k+1}= \cTr(\cK_{\theta,k})$, $\Or^k$, the refinement region in $\cT_k$, and denote by $u_k$ the solutions 
of \eref{normalT} in $\U_k=\U_{\cT_k}$. Now consider $m$ refinement cycles, giving 
$$
u_{k+j}, \,\, \cT_{k+j},\,\, w_{k+j}= u_{k+j+1}-u_{k+j}, \,\, \cR_{k+j}= f- \cB u_{k+j},\quad j=0,\ldots m-1.
$$
By orthogonality
\begin{align}
\label{orthog}
&\|u^*- u_{k+m}\|_\U^2 = \|u^*- u_k\|^2_\U- \|w_{k+m-1}\|^2_\U - \cdots - \|w_k\|^2_\U   .
\end{align}
Fix   $\bar\rho<1$ and suppose that for some $m>l\ge 1$ one has $\|\cR_{\cT_{k+m-l}}\|_\Om
\le   \bar\rho \|\cR_{\cT_{k+m-l-1}}\|_\Om$.
Then
 $$
 \begin{aligned}
 \|u^*-u_{k+m}\|_\U&\le \|u^*- u_{k+m-l}\|_\U=  \|\cR_{\cT_{k+m-l}}\|_\Om\le \bar\rho 
 \|\cR_{\cT_{k+m-l-1}}\|_\Om\\
 & =\bar\rho\|u^*- u_{k+m-l-1}\|_\U\le \bar\rho  \|u^*- u_{k }\|_\U,
 \end{aligned}
 $$
which is the claimed error reduction by a factor $\bar\rho<1$.
Suppose therefore that
\be
\label{suppose}
\|\cR_{\cT_{k+m-l}}\|_\Om
>   \bar\rho \|\cR_{\cT_{k+m-l-1}}\|_\Om,\quad \forall\,  1\le l <m.
\ee
From \eref{osc4} and the above observations we know that with $\W_{k+m-j}:= \W(\cT_{k+m-j};\cT_{k+m+1-j})$
\begin{align}
\label{mdiff}
\sum_{j=1}^{m}\|w_{k+m-j}\|^2_\U&\ge \sum_{j=1}^{m}\xi^2 \|\cR_{k+m-j}\|^2_\Om- \osc_{\cT_{k+m-j}}(f,\W_{k+m-j})^2
\nonumber\\
&\ge \xi^2 (\bar\rho^2 +\cdots +\bar\rho^{2(m-1)})\|\cR_{\cT_k}\|^2_\Om - \sum_{j=1}^{m}
\osc_{\cT_{k+m-j}}(f,\W_{k+m-j})^2\nonumber\\
&\ge \xi^2 (\bar\rho^2 +\cdots +\bar\rho^{2(m-1)})\|\cR_{\cT_k}\|^2_\Om - 
C_0 \|\cR_{\cT_k}\|^2_\Om\nonumber\\
&= \Big\{\frac{\bar\rho^2\xi^2(1-\bar\rho^{2(m-1)})}{1-\bar\rho^2}-C_0\Big\}\|\cR_{\cT_k}\|^2_\Om .
\end{align}
Thus, whenever   $m,\bar\rho$ satisfy  for some $\eta>0$  
\be
\label{C0eta}
\frac{\bar\rho^2\xi^2(1-\bar\rho^{2(m-1)})}{1-\bar\rho^2}\ge C_0+\eta.
\ee
We infer from \eref{orthog},  \eref{mdiff}, and Remark \ref{rem:iff} that
$$
\|u^*-u_{k+m}\|_\U\le \sqrt{1-\eta}\|u^*-u_k\|_\U.
$$
On account of the above alternative, we have
$$
\|u^*-u_{k+m}\|_\U\le \max\big\{\bar\rho,\sqrt{1-\eta}\big\}\|u^*-u_k\|_\U.
$$
Now, choosing $\eta=1-\bar\rho^2$, \eref{C0eta} takes the form
$$
 \bar\rho^2\Big\{\frac{\xi^2(1- \bar\rho^{2(m-1)})}{1- \bar\rho^2}+1\Big\}\ge C_0+1.
$$
By Lemma \ref{lem:xiC0}, there indeed exist  $\rho=\bar\rho(m)<1,m\in\N$ such that this inequality
holds,
which completes the proof.\hfill $\Box$\\

\section{Towards a Practical Scheme}\label{sec:practical}
%%%%%%%%%%%%%%%%%%%%%%%%%%%%%%%%%%%%%%%%%%%%%%%%%%%

The  obvious remaining essential question is for which $f$ does \eref{ell2} hold?
This is intimately tied to the perhaps  main idealization in the previously discussed scheme, namely 
 the
assumption of having exact access to $f$ and resulting quantities like residuals. These quantities
are subject to quadrature errors. We re-interpret the entailed perturbation of quantities, involving $f$, as resulting from an exact integration with $f$ replaced by a piecewise polynomial approximation.
 We therefore also assume 
in what follows that the coefficients in $\cB$ are piecewise polynomial on the current partition.
We suppress for simplicity of exposition a correspondingly extended perturbation analysis, taking
piecewise polynomial approximations of such coefficients into account, e.g. absorbed by inexact solves
that can be handled along the lines in \cite{robfocm}.
 
So, let us proceed under the
  working assumption  that
\be
\label{working}
(\t f,\cB v)_\Om, \quad\|\cR_\cT(v,\t f)\|_{\Om'},\quad \cR_\cT(v,\t f):= \t f - \cB v, \quad v\in \U_\cT,
\ee
can be computed exactly and that $\textsc{Solve}[\t f,\cT]\to \t u_\cT$ can be carried out exactly
when $\t f$ is a piecewise polynomial.
%%%%%%%%%%%%%%%%%%%%%%%%%%%%%%%%%
\subsection{Piecewise Polynomial Loads and a Reduction Block}
%%%%%%%%%%%%%%%%%%%%%%%%%%%%%%%%%%% 

Computing a piecewise polynomial approximation to a given $f\in \L$ would, of course also
require the ability to evaluate $f$ in some way, e.g. to employ sufficidently accurate quadrature. Proposing a concrete procedure depends again on the specific case at hand.
We will therefore work here under the following specified assumption analogous to \cite{robfocm}.\\

\noindent
Given $f\in \L$ and a partition $\cT\in\mT$ we have at hand a routine\\[2mm]
$\textsc{Appr}[f,\cT,\delta]\to ( f_\delta,\cT^f_\delta)$ \\[2mm]
where $\cT^f_\delta\in \mT$ is a refinement of $\cT$, and $f_\delta$ is a piecewise
polynomial on $\cT^f_\delta$  such that
\be
\label{fdelta}
\|f- f_\delta\|_\Om \le \delta.
\ee
Here $\cT^f_\delta$ is obtained from $\cT$ through a finite number of calls of $\textsc{Subd}[\cT,\cdot,1]$, followed by a closing step to restore $\mT$-compliance.

For a given piecewise polynomial load, it will be useful to view  a certain fixed number of 
successive cycles in
\textsc{Alsm-Ideal} as a single {\em refinement block} which we formulate as the following sub-routine.
 \\

\noindent
$\textsc{Reduce}[\cT,\t f,\t u,\theta, \rho]\to  \cTr, \bar u_\cTr\in \U_{ \cTr}$\\[2mm] 
Input: $\cT\in \mT$, $\t f$ piecewise polynomial on $\cT$, $\t u= \textsc{Solve}[\t f,\cT]$, $\theta\in (0,1)$, $\rho<1$;
\begin{itemize}
\item[{~}] 
Set  $\t\cT=\cT$;\\
do:\\
\hspace*{3mm} $\textsc{Mark}[\t \cT,\t u,\t f,\theta]\to \cK$; \\ 
\hspace*{3mm} $\textsc{Refine}[\t\cT,\cK]\to \cT$;\\ 
\hspace*{3mm} $\textsc{Solve}[\t f,\cT]\to \bar u$;
\\
\hspace*{1cm}while
\be
\label{while}
\|\t f- \cB\bar u \|_\Om > \rho\|\t f-\cB \t u\|_\Om,
\ee
\hspace*{1cm} repeat;\\
\hspace*{3mm} else, output $ \cTr=\cT$, $\bar u_\cTr =\bar u $.
\end{itemize}
\bigskip

By construction, whenever \textsc{Reduce} terminates, one has   
\be
\label{reduce}
\|\t f - \cB  \bar u_\cTr\|_\Om \le \rho\|\t f-\cB\t u\|_\Om.
\ee

It will be shown later below that when $f$ is piecewise polynomial on a current partition,
\textsc{Reduce} indeed terminates   within a fixed finite number $m$   of inner calls of \textsc{Refine}
where $m$ satisfies \eref{solution}.

%%%%%%%%%%%%%%%%%%%%%%%%%%%%%%%%%%%%
\subsection{\textsc{ALSM} and the Main Result}
%%%%%%%%%%%%%%%%%%%%%%%%%%%%%%%%

To formulate the main adaptive scheme let $m\in\N$, $\rho<1$ satisfy \eref{solution} and    choose $\omega\in (0,1)$
such that
\be
\label{omega}
  \hat\rho := \frac{\omega + \rho}{1-\omega}  <1.
\ee
\begin{rem}
\label{rem:theta}
We have to keep in mind that $\rho$ from \eref{solution}, hence $\hat\rho$, depend on $\theta$ (through $\xi(\theta)$) but exist for any $\theta<1$. If this dependence matters, we write $\rho=\rho(\theta), \hat\rho=\hat\rho(\theta)$.
\end{rem}

The idea behind the following scheme is: first generate a piecewise polynomial approximation 
to the load $f$ by refining the current partition   until the approximation error is dominated
by the residual of the Galerkin solution on the updated partition: then apply a reduction block
\textsc{Reduce} which is then guaranteed to terminate after at most $m$ calls of \textsc{Refine}.\\

\noindent
$\textsc{Alsm}[\cB,f,\theta, \e] \to \cT_\e,\, u_\e\in \U_{\cT_\e}$:

 \begin{itemize}
\item[(A0)]
Set $\e>0$;\\
choose $\theta< 1$, $\rho$, $\omega$ according to \eref{solution},\eref{omega},
respectively;\\
choose $\cT_0$, $u_0=0$,  $\ell\in\N$,   
$k=0$;\\ 
fix $\delta_0 = \| f\|_\Om$,   $\delta_0\to \delta$;\\
\item[(A1)]
set $\cT_k\to \cT$, $\delta_k\to \delta$;\\
$\textsc{Appr}[f,\cT,\delta ]\to [  f_{\delta },  \t\cT_k]$;\\
 $\textsc{Solve}[  f_{\delta }, \t\cT_k ]\to \bar u_{\delta }$; \\[2mm]
while 
\be
 \label{delom}
 \delta  +  \|f_{\delta }- \cB \bar u_{\delta }\|_\Om >\e,
\ee
\hspace*{2cm}
while
\be
\label{while2}\delta >\omega\| f_{\delta } - \cB\bar u_{\delta }\|_\Om;
\ee
\hspace*{3cm} $\delta /2\to \delta_k$, $\t\cT_k\to \cT_k$,  and go to (A1);\\ 
\hspace*{2.5cm}
\hspace*{-5mm}$f_{k}=  f_\delta, \bar u_{k} =\bar u_\delta$, $\delta_k=\delta$, go to (A2);\\[2mm]
stop and set $\cT_\e=  \t\cT_{k}$, $u_\e= \bar u_\delta$;

 \item[(A2)]
   $\textsc{Reduce}[\t\cT_{k } ,f_{ k  },\bar u_k, \theta,\rho ]\to \cT_{k+1}, u_{k+1}$;\\
 while $\delta_k +  \|f_k - \cB   u_{k+1}\|_\Om > \e$;\\
 \hspace*{2cm} 
  set $k+1\to k$ and goto (A1);\\[2mm]
  stop and accept
 $\cT_\e = \cT_{k+1}$, $u_\e =  u_{k+1}$ as the output.
 \end{itemize}
 
 The analysis of the above scheme hinges on the validity of the following standard approximation properties of the 
trial spaces, see \cite[Section 2.5.2]{BBF} for known relevant scenarios.\\
 
%%%%%%%%%%%%%%%%%%%%%%%%
\newcommand{\CAP}{C_{\mbox{\tiny{CAP}}}} 

\noindent
{\bf Assumption CAP:} (Canonical Approximation Property) {\it For  $\Or\subset \Om$ as above
we assume that the trial spaces enjoy canonical approximation properties in the sense that the following
localized Bramble-Hilbert Lemma holds: if $f\in H^s(\Om;\R^q)$ for some $s>0$,
there exists a constant $\CAP$, independent of $\Or$, $f$, such that 
\be
\label{CAP}
\inf_{v\in \U(\Or;\cTr)}\|\cB(u^*_{\Om'}-v)\|_{\Or}
\le \CAP \sum_{K\in\cTr,K\subset \Or} h_K^{s}\|f\|_{H^s(K;\R^q)}.
\ee
} 

\begin{theorem}
\label{thm:pract1}
Under Assumption (CAP), one has:\\[2mm]
(a)  The output of \textsc{Alsm}  meets the target tolerance $\|u^*- u_\e\|_\U\le \e$.\\[2mm]
(b) \textsc{Reduce} terminates in each of its calls  within at most $m$ inner calls of \textsc{Refine},
where $m, \rho$ are related by \eref{solution}.\\[2mm]
(c) For $\hat\rho<1$ as in \eref{omega}  
the sequence $u_k,\cT_k$ produced by $\textsc{Alsm}[\cB,f,\theta,\e]$ satisfies,      as long as $u_{k+1}\neq u_\e$,     
\be
\label{link}
\|u^*- u_{k+1}\|_\U\le \hat\rho\|u^*-u_k\|_\U.
\ee
Hence, the scheme terminates  after $O(|\log \e|)$
number of calls of \textsc{Reduce}.
 \end{theorem}

 Note that (A1) is the only step in \textsc{Alsm} where the termination of the inner loop is not quantified.  In order to branch into (A2) the piecewise polynomial approximation error to $f$ must be smaller than a certain fraction of the current residual. Without any
 further assumption of $f$, it  might happen that \eref{delom} fails before \eref{while2} does,
 in which case the scheme also terminates and, as shown below, also meets the target tolerance. 
 Thus, the theorem asserts a linear error reduction with regard to the calls of \textsc{Reduce}
 but so far cannot quantify the refinements required to resolve $f$ sufficiently well by piecewise
 polynomial approximations.  

%%%%%%%%%%%%%%%%%%%%%%%%%%%%%%%%%%%%%%%%%%
\subsection{Proof of Theorem \ref{thm:pract1}:}
%%%%%%%%%%%%%%%%%%%%%%%%%%%%%%%%%%%%%%%%%%%%%%%%

Our first step is to show   that
when the input load is a piecewise polynomial on $\cT$ the routine \textsc{Reduce} indeed terminates 
for an appropriate $\rho$ after a uniformly bounded finite number of steps, provided that
the trial spaces $\U_\cT$, $\cT\in\mT$ have a canonical approximation property.

 We make frequent use of the following:
\begin{remark}
\label{rem:proju*}
For any $u\in\U$ the unique solution $u_\Or\in \U(\Or)$ of
\be
\label{locprob}
(\cB u_\Or,\cB v)_\Or = (\cB u,\cB v)_\Or,\quad v\in \U(\Or).
\ee
is given by 
$P_{\U(\Or)}u\in \U(\Or)$.  
Hence
\be
\label{Lloc}
\cB (P_{\U(\Or)}u) = \cB u_\Or = (\cB u)|_\Or \quad\mbox{in $\,\,\L(\Or)$}.
\ee
\end{remark}
\noindent
{\bf Proof:} By definition,  
\be
\label{Urdef}
(P_{\U(\Or)}u, v)_\U = (u,v)_\U \quad \forall\, v\in \U(\Or).
\ee
The left hand side of \eref{Urdef} agrees, by definition, with
$$
(\cB P_{\U(\Or)}u, \cB v)_\Om = (\cB P_{\U(\Or)}u, \cB v)_\Or,\quad \mbox{ (since $v\in \U(\Or)$)}
$$
while the right hand side of \eref{Urdef} can be rewritten as
$$
(u,v)_\U= (\cB u,\cB v)_\Om= (\cB u,\cB v)_\Or,
$$
again since $v\in \U(\Or)$. By well-posedness of \eref{locprob}, we must have $u_\Or=P_{\U(\Or)}u$,
as claimed. \hfill $\Box$\\

We show next that   $\osc_\cT(f,\W)$ can be bounded by a best approximation error.
\begin{lemma}
\label{lem:ba}
Adhering to the above notation, we have 
\be
\label{oscest}
\osc_\cT(f,\W)\le  \inf_{\t u\in \U(\Or; \cTr)}\|P_{\U(\Or)}(u^*-v_\cT)- \t u\|_{\U(\Or)}\quad \forall\,
v_\cT\in\U_\cT.
\ee
\end{lemma}
\noindent
{\bf Proof:}
 We infer from the definitions
of $\Wi$ and $\W$, \eref{Wdef}, \eref{Wi}, that
$$
\begin{array}{c}
P_{\Wi} u= (I- P_{\U_\cT})P_{\U_\cT(\Or)}u,\\ P_{\W} u= (I- P_{\U_\cT})P_{\U_{\cTr}}u
= (I- P_{\U_\cT})P_{\U_\cT(\Or; \cTr )}u  ,
\end{array}
$$
where we have used \eref{write} in the last step.
Hence,
\be
\label{diffW}
\begin{aligned}
(P_{\Wi}-P_\W)u&= (I- P_{\U_\cT})(P_{\U_\cT(\Or)}- P_{\U_\cT(\Or;\cTr)})u\\
& = (I- P_{\U_\cT})(P_{\U (\Or)}- P_{\U (\Or;\cTr)})u\quad u\in\U.
\end{aligned}
\ee
Now, recall from Lemma \ref{lem:rep} that
$\osc_\cT(f,\W)= \|(P_{\Wi}-P_\W)u^*\|_\U$.  
Noticing that $(P_{\Wi}-P_\W)u^*= (P_{\Wi}-P_\W)(u^*-v_\cT)$ for any $v_\cT\in\U_\cT$ and     invoking \eqref{diffW},  provides
$$
\begin{aligned}
\osc_\cT(f,\W)^2&=\|(P_{\Wi}-P_\W)u^*\|^2_\U\\
&   = ((I- P_{\U_\cT})(P_{\U(\Or)} 
  - P_{\U(\Or;\cTr)})(u^*-v_\cT)  ,(P_{\Wi}-P_\W)u^*)_\U\\
&\qquad  \le \|(P_{\U(\Or)} -P_{\U(\Or;\cTr)})(u^*-v_\cT) \|_{\U(\Or)} \|(P_{\Wi}-P_\W)u^*\|_\U\\
&= \inf_{\t u\in \U(\Or;\cTr)}\|P_{\U(\Or)} (u^*-v_\cT)- \t u\|_{\U(\Or)}\|(P_{\Wi}-P_\W)u^*\|_\U,
\end{aligned}
$$ 
whence the assertion follows.\hfill $\Box$\\

We adhere to the previous notation and denote
by $\cT=\cT_k$, $\cT_{k+1}= \cTr(\cK_{\theta,k})$, $\Or^k$, the refinements and approximate solutions,
generated by \textsc{Reduce} (until termination). Accordingly,   $u_k$ stands for the solution 
of \eref{normalT} in $\U_k=\U_{\cT_k}$. Moreover,  we abbreviate $\W_{k+m-j}:= \W(\cT_{k+m-j},\cT_{k+m+1-j})$,
$\osc_{k+m-j}(f):=\osc_{\cT_{m+k-j}}(f,\cT_{m+k-j})$, $j=1,\ldots,m$.

\begin{lemma}
\label{lem:CAP}
 If {\bf Assumption CAP} holds and $f$ is a piecewise polynomial on $\cT=\cT_k$. then for any $m\in\N$ 
\be
\label{oscsum}
\sum_{j=1}^m \osc_{m+k-j}(f)^2 \le C_0\|\cR_{\cT_k}\|^2_\Om,
\ee
where the constant $C_0$ depends only on $C_{inv},\CAP,\kappa,s,\mT,\Om$
(see {\bf Assumption CAP}). Here   $C_{inv}$
is the constant appearing in a local inverse estimate for the trial spaces, see \cite{5authors}.
 \end{lemma}
\noindent
{\bf Proof:}  First fix $\cT$ and set 
   $e^*_\cT= u^*-u_\cT$, so that  $\cB e^*_\cT= \cR_{\cT_k}$ and we can apply  \eref{oscest} to
  $u^*-v_\cT= e^*_{\cT}$. In view of  
 Remark \ref{rem:proju*}, this yields
\be
\label{oscest2}
\osc_{\cT}(f,\W) \le  \inf_{\t u\in \U (\Or;\cTr )}\|\cR_\cT - \cB \t u\|_{\Or }.
\ee
We will apply this to $\cT=\cT_{m+k-j}, \Or^{m+k-j}$, $j=1,\ldots,m$.
Now suppose that $f$ is a piecewise polynomial subordinate to the partition $\cT_k $. Hence $\cR_{\cT_k}$ is a 
piecewise polynomial on $\cT_k$ as well as on the refinements $\cT_{k+m-j}$, $j=1,\ldots,m$.
Even, when $\cR_{\cT_k}$ is discontinuous there exists a positive $s$ so that $\cR_{\cT_k}\in H^s(\Omega;\R^q)$. Note that the   cells in $\Or^{k+m-j} $ arise for $j\le m$   from at least one subdivision of a cell in $\cT_k$.

To proceed we shall make frequent use of the weak sub-additivity of
the Sobolev semi-norms in the form
\be
\label{subadd}
\sum_{K\in \cK}|u|^2_{H^s(K;\R^q)}\le  |u|^2_{H^s(\Om_\cK)},
\ee
where $\cK\subset \cT$ and $\Om_\cK=\bigcup\,\{K: K\in \cK\}$. In view of 
\eqref{subadd}, we have
\be
\label{merge}
\sum_{\substack{K'\subset K\\ K' \mbox{ children of }\,K}}h_{K'}^{2s}|\cR_{\cT}|^2_{H^s(K')}
\le   \kappa^{2s}h_K^{2s}|\cR_{\cT}|^2_{H^s(K)},\quad K\in\cT.
\ee

Consider now the collection of cells $K'\in \cT_{m+k+1-j}$ that are contained in 
$\Or^{m+k-j}$. Then, such a $K'$ has a unique {\em parent cell} $K=K(K')\in \cT_{k+m-j}$. Denoting 
by $h_{K'}:={\rm diam}\,K'$ the local meshsize of $K'$, we have by assumption $h_{K'}\le \kappa h_K$
for a fixed $\kappa\in (0,1)$, depending only on $\mT$.

Now we invoke {\bf Assumption CAP} to conclude that 
\begin{align}
\label{oscest3}
 \osc_{m+k-j}(f)^2 
&   \le
  \inf_{\t u\in \U(\Or^{k+m-j};\cT_{m+k+1-j})}\|\cR_{\cT_k} - \cB \t u\|^2_{\Or^{m+k-j}}\nonumber \\
&   \le \CAP \sum_{K'\in\cT_{m+k+1-j},K'\subset\Or^{m+k-j}}h_{K'}^{2s} |\cR_\cT|^2_{H^s(K';\R^q)}.
\end{align}
Hence
\begin{align}
\label{sumosc}
 \sum_{j=1}^{m}
\osc_{m+k-j}(f)^2 &\le \CAP   \sum_{j=1}^m  \sum_{K\in\cT_{K+m+1-j},K\subset\Or^{k+m-j}}h^{2s}_K |\cR_\cT|^2_{H^s(K;\R^q)}.
\end{align}
We proceed with partitioning the total refinement region
$$
\wh\Or^k:= \bigcup_{j=1}^m\Or^{m+k-j},
$$
into the collection of cells
$$
\cK_l:=\{K\in \cT_{k+l}: K \subset\Or^{k+l},\,K\cap \Or^{k+l'}=\emptyset, l'<l\},\quad l=0,\ldots, m-1,
$$
so that obviously 
\be
\label{Omdeco}
\wh\Or^k= \bigcup_{l=0}^{m-1}\Om(\cK_l),\quad \Om(\cK_l)\cap \Om(\cK_q)=\emptyset ,\quad q\neq l,
\ee
where we denote as before for any collection $\cK$ of cells $\Om(\cK):=\bigcup\{K\in\cK\}$. 
 Note next that when $K\in \cK_l$ intersects any of the regions $\Or^{ k+l'}$, for some $l'>l$
one must have $K\subset\Or^{ k+l'}$. Define for $K\in \cK_l$
$$
d_{K,l}:= \#\{ l'>l: K\subset \Or^{k+l'}\}
$$
the number of times $K$ is covered by a higher level refinement region. Hence $d_{K,l}$ is for every
$K\in \cK_l$, $l=0,\ldots,m-1$, the number of times $K$ got subdivided in the process.
Therefore, invoking \eref{subadd} and \eref{merge}, we obtain
\be
\label{oscestfinal}
 \sum_{j=1}^{m}
\osc_{m+k-j}(f)^2 \le \CAP \sum_{l=0}^{m-1} \sum_{K\in \cK_l}\Big(\sum_{r=0}^{d_{K,l}}\kappa^{2s r}\Big)h_K^{2s}|\cR_{\cT_k}|^2_{H^s(K;\R^q)}.
\ee
Since $\kappa<1$ and since, by \eref{Omdeco}, the $\cK_l$ mark disjoint regions, we obtain, after
using an inverse estimate
\be
\label{oscestfinal+}
 \sum_{j=1}^{m}
\osc_{m+k-j}(f)^2 \le \CAP C_\kappa C_{inv} \sum_{l=0}^{m-1} \|\cR_{\cT_k}\|^2_{\Om(\cK_l)}
=: C_0 \|\cR_{\cT_k}\|^2_{\wh\Or^k},
\ee
which completes the proof of Lemma \ref{lem:CAP}. \hfill $\Box$\\
 
 We are now ready to complete the proof of Theorem \ref{thm:pract1}.\\

\noindent
{\bf Ad (a):} If in (A1), due to failure of \eref{delom}, we have $\delta_k+ \|f_{ k}- \cB \bar u_{\delta}\|_\Om \le \e$ we obtain
$$
\|u^*-\bar u_{\delta}\|_\U\le \|f- f_{ k}\|_\Om + \|f_{ k}-\cB \bar u_{\delta}\|_\Om \le \delta_k + \|f_{ k}-\cB \bar u_{\delta}\|_\Om
\le \e,
$$
so that, by definition of $u_\e$ in this case, the scheme meets the target accuracy upon termination.
 \\[2mm]
{\bf Ad (b):} Since the load in \textsc{Reduce} is a piecewise polynomianl on the input partition, Lemma \ref{lem:CAP} says that \eref{ell2} holds. The claim follows then  
from Theorem \ref{thm:main1} 
  with $m,\rho$ satisfying \eref{solution}.\\[2mm]
{\bf Ad (c):} Suppose the scheme terminates in (A2). The input $\bar u_k, f_k$ of \textsc{Reduce} in (A2)
must have failed \eqref{while2}, hence satisfies $\delta_k\le \omega\|f_k- \cB \bar u_k\|_\Om$.
Under Assumption (CAP), if $u_{k+1}, \cT_{k+1}$ are the output of \textsc{Reduce} in step (A2),
we know from  Theorem \ref{thm:main1} together with Lemma \ref{lem:CAP} 
 that the exact solution $u^*_{k }:= \cB^{-1}f_{ k }$ satisfies $\|u^*_{k }-u_{k+1}\|_\U \le \rho \|u^*_{k }-\bar u_k\|_\U$, where, as a result of (A1),
  $\textsc{Solve}[f_{k},\cT_{k}]=\bar u_{ k }$. Hence
\begin{align*}
\|u^*- u_{k+1}\|_\U &\le \|u^*- u^*_{k }\|_\U+ \| u^*_{k }- u_{k+1}\|_\U 
 = \|f- f_{ k }\|_\Om + \| u^*_{k }- u_{k+1}\|_\U\\
&\le \delta_k + \rho \|u^*_{k }-\bar u_k\|_\U 
 \stackrel{\eref{while2}}{\le} \omega \|f_{{k }} - \cB \bar u_k\|_\Om + \rho \|u^*_{k }-\bar u_k\|_\U\\
&= (\rho +\omega)\|u^*_{k }-\bar u_k\|_\U.
\end{align*}
On the other hand, 
\begin{align*}
\|u^*_{k }-\bar u_k\|_\U &\le  \|u^*_{k }- u^*\|_\U + \|u^*-\bar u_k\|_\U 
 = \|f_{ k }- f\|_\Om + \|u^*-\bar u_k\|_\U\\
&\le \delta_k + \|u^*-\bar u_k\|_\U 
 \le \omega\|u^*_{k }-\bar u_k\|_\U + \|u^*-\bar u_k\|_\U,
\end{align*}
which, upon combining with the preceding bound, yields
\be
\label{zack}
\|u^*- u_{k+1}\|_\U\le \frac{\omega+\rho}{1-\omega}\|u^*-u_k\|_\U .
\ee
In view of \eref{omega}, this confirms \eref{link} and finishes the proof of (c).
\hfill$\Box$\\

 Theorem \ref{thm:pract1} does still not quantify the total work yet, needed to meet the target accuracy. It cannot do so because
 we do not know when the while-loop \eref{while2} terminates, i.e., how fast $f_\delta\to f$ and how fast $\|f- \cB \t u_\delta\|_\Om\to 0$ in step (A1). Such quantifications require additional information 
 about $f$ 
  beyond just $f\in \L$, see Remark
 \ref{rem:f'}.
In the next section we discuss exemplary conditions on $f$ that lead to complexity bounds.

%%%%%%%%%%%%%%%%%%%%%%%%%%%%%%%%%%%%%%%%%%%%%
\section{Complexity - a Light Version}\label{sec:compl}
%%%%%%%%%%%%%%%%%%%%%%%%%%%%%%%%%%%%
We continue focusing solely on piecewise polynomial perturbations of $f$ because  this suffices first, to
bring out the role of
 $\osc_\cT(f,\W)$, and  second, to  highlight those aspects where arguments have to differ from
 the lines in \textsc{Afem}.
In principle, further perturbations, caused e.g. by more general coefficients in $\cB$
for which $\cB u_\cT$ is no longer piecewise polynomial, or inexact solves can be treated
in a similar fashion, following essentially the strategy in \cite{robfocm}. To account for
these further perturbation in complexity estimates would require - on the current level of generality -
corresponding assumptions that need to be checked in each concrete case at hand.
 
Again, since we haven't specified the structure of $\mT$, we work under the following additional 
requirement which is known to hold in many relevant scenarios and plays the same role in classical 
\textsc{Afem} scenarios, see e.g. \cite{BDD,S}.\\[2mm]
 {\bf RSP:} (Refinement Summation Property) {\it
One has for $\cK_j=\cK(\cT_j,\theta)$
\be
\label{RSP}
\#\cT_N-\#\cT_0\lsim \sum_{j=0}^{N-1}\#\cK_j,
\ee 
where the proportionality constant is independent of $N$.}\\

Finally, we introduce the {\em approximation classes} $\cA^r= \cA^r_\L\times \cA^r_\U$, for some $r>0$, comprised of
all pairs $(f, \cB^{-1}f)\in \L\times \U$ with the following property: \\[2mm]
(P1): there exists
a $C_f<\infty$ such that
\be
\label{Af}
N^r \inf_{\substack{\cT\in\mT\\(\#\cT -\#\cT_0)\le N}}\|f- f_\cT\|_\L\le C_f,
\ee
where $f_\cT$ is the $\L$-projection of $f$ to the space of all piecewise 
polynomials $\PP_\cT$ in $H^s(\Omega;\R^q)$ subordinate to $\cT\in\mT$.\\[2mm]
(P2)
For $u=\cB^{-1}f$, there exits a constant $C_u<\infty$ such that
\be
\label{Au}
N^r  \inf_{\substack{\cT\in\mT\\(\#\cT -\#\cT_0)\le N}}\|u- u_\cT\|_\U\le C_u,
\ee

The smallest such constants $C_f, C_u$ for which \eref{Af}, \eref{Au} hold respectively,
define semi-quasi-norms $|f|_{\cA^r(\L)}$, $|u|_{\cA^r(\U)}$, so that $\|f\|_\L+|f|_{\cA^r(\L)}$ and $\|u\|_\U +|u|_{\cA^r(\U)}$ are (quasi-) norms on $\L, \U$, respectively.

Hence, given $\e>0$, there exists a $\cT_\e\in\mT$ such that $f\in \cA^r(\L)$ can be approximated in $\L$ by a piecewise polynomial $\PP_{\cT_\e}$ to accuracy $\e$, i.e.,  $\|f- f_{\cT_\e}\|_\L\le \e$, 
 where for some multiplicatice constant independent of $f,$
\be
\label{fT}
\#\cT_\e-\#\cT_0\lsim |f|^{1/r}_{\cA^r(\L)}\e^{-1/r}.
\ee
An  analogous relation holds for $u$.\\

We follow \cite{robfocm} and assume that the routine $\textsc{Appr}$ has the following
property which can only hold when $f\in \cA^r$.

\medskip
\noindent
{\bf Assumption Af:} {\it The pair $(f,\textsc{Appr})$ is $r$-{\em optimal},
i.e., for a given $\cT$, the  output $(\hat\cT ,f_{\hat\cT})$ of $\textsc{Appr}[f,\cT  ;\delta]$
satisfies for any $\delta>0$
\be
\label{foptimal}
\#\hat\cT - \#\cT \le c(f)^{1/r} \delta^{-1/r},\quad \|f- f_{\hat\cT}\|_\Om\le \delta.
\ee
where $c(f)$ depends only on $f$.
}\\
 
 For instance, if suitable sufficiently accurate quadrature methods are available
 one can employ
the Binev-DeVore tree algoritm (in a perturbed form) to generate $f_\delta$, see e.g.
\cite{CP}.

 \begin{theorem}
\label{thm:pract2}
Adhering to the above notation, the following hold:\\[2mm]
(i) Under Assumptions {\bf  Af}, and {\bf CAP}, Algorithm $\textsc{Alsm}[\cB,f,\theta,\e]$ terminates after finitely
many steps, and outputs a partition $\cT_\e\in\mT$ such that
\be
\label{target}
\|u^*- u_{\cT_\e}\|_\U \le \e.
\ee
\medskip
(ii) Suppose in addition that Assumption {\bf RSP} holds and that $(f,u^*=\cB^{-1}f)\in \cA^r:=\cA^r(\L)\times 
\cA^r(\U)$. Then, for any $\theta\in (0,1)$ 
  the final partition $\cT_\e$, produced by $\textsc{Alsm}[\cB,f, \theta, \e]$, satisfies
\be
\label{Tsize}
\#\cT_\e - \#\cT_0 \lsim \e^{-1/r}|(f,u^*)|^{1/r}_{\cA^r},\quad \e >0,
\ee
with a constant   depending on $\theta,\mT$ and on the constants in Assumptions {\bf Af, CAP, RSP}.
\end{theorem}
\begin{remark}
\label{rem:theta}
There is no constraint on $\theta$ other than $\theta\in (0,1)$.
\end{remark}
%%%%%%%%%%%%%%%%%%%%%%%%%%%%%%%%%%%%%
\subsection{Proof of Theorem \ref{thm:pract2}:}
%%%%%%%%%%%%%%%%%%%%%%%%%%%%%%%%%%%%%%%%%%%%%
The proof of Theorem \ref{thm:pract2} is  based
on three Lemmata.
The first one concerns the role of D\"orfler marking in the present context.

\begin{lemma}
\label{lem:reduce}
 Assume that $\t\cT\succ \cT$ is any refinement with refinement region $\t\Om_{\rm ref}$ such that
\be
\label{tilde}
\|u^*- u_{\t\cT}\|^2_\U \le \lambda\|u^*-u_\cT\|^2_\U\quad \mbox{holds for some }\,\,
\lambda < 1  .
\ee  
 Then, one must have
\be
\label{tbulk}
\|\cR_\cT\|^2_{\t\Om_{\rm ref}}\ge (1-\lambda) \|\cR_\cT\|^2_\Om,
\ee
in other words the refinement must obey a  bulk criterion.\\[2mm]
\noindent
\end{lemma}
\noindent
{\bf Proof:}
Let $\t\Om^c_{\rm ref}:=\Om\setminus \t\Om_{\rm ref}$, $\t\cR^c_\cT:= \chi_{\Om\setminus \t\Om_{\rm ref}}\cR_\cT$. Consider  the variational problem: find $z_\cT\in\U$ such that
\be
\label{zT}
(\cB z_\cT,\cB v)_\Om= (\t\cR_\cT^c,\cB v)_\Om,\quad v\in \U.
\ee
This is a well posed problem and
\be
\label{znorm}
\|z_\cT\|_\U=\|\cB z_\cT\|_\Om = \|\t\cR_\cT^c\|_\Om =\|\cR_\cT\|_{\t\Om^c_{\rm ref}}.
\ee
By definition,
\begin{align*}
\|u^*- u_{\t\cT}\|_\U &= \sup_{v\in \U}\frac{((I- P_{\t\cT})u^*,  v)_\U}{\|v\|_\U} =
\sup_{v\in \U}\frac{((I- P_{\t\cT})(u^*-u_\cT),  v)_\U}{\|v\|_\U} \\
&= \sup_{v\in \U}\frac{(\cR_\cT, \cB(I-  P_{\t\cT}) v)_\Om}{\|\cB v\|_\Om}\\
 &\ge \frac{(\cR_\cT,\cB((I- P_{\t\cT})z_\cT)_\Om}{\|\t\cR_\cT^c\|_\Om}.
\end{align*}
Now, since $\cT \prec \t\cT$, hence $\U_\cT\subset \U_{\t\cT}$
\begin{align*}
 (\cR_\cT,\cB((I- P_{\t\cT})z_\cT)_\Om
 &= (z_\cT,(I- P_{\U_{\t\cT}})(u^*-u_\cT))_\U = (z_\cT,u^*- u_\cT)_\U\\
&= (\cB z_\cT, \cB(u^*-u_\cT))_\Om= (\t\cR_\cT^c,\cB(u^*-u_\cT))_\Om\\
&= (\t\cR_\cT^c,\cR_\cT)_\Om = \|\cR_\cT \|^2_{\Om\setminus\t\Om_{\rm ref}},
\end{align*}
which gives
$$
\|u^*- u_{\t\cT}\|_\U\ge  \|\cR_\cT \|_{\Om\setminus\t\Om_{\rm ref}}.
$$
Now suppose that   $\|\cR_\cT\|^2_{\t\Om_{\rm ref}}=\kappa
\|\cR_\cT\|^2_\Om$ for some $\kappa <1$. Then $(1-\kappa)\|\cR_\cT\|^2_\Om =  
\|\cR_\cT \|^2_{\Om\setminus\t\Om_{\rm ref}}$. Thus,
$$
\|u^*- u_{\t\cT}\|^2_\U\ge (1-\kappa)\|u^*-u_\cT\|^2_\U ,
$$
and \eref{tilde} says that $\lambda \ge (1-\kappa)$ which gives \eref{tbulk}.
 \hfill $\Box$\\
 
The next lemma concerns the complexity of step (A1) in \textsc{Alsm} when $f$ satisfies Assumption {\bf Af}.
\begin{lemma}
\label{lem:RHS}
If $(f,\textsc{Appr})$ is $r$-optimal
the output 
$
\cT_k , \, f_{ k} ,\,  \bar u_k 
$ of step (A1) in \textsc{Alsm}  
satisfies
\be
\label{outRHS}
\#\t\cT_k-\#\cT_k \le C  c(f)^{ 1/r} \delta_k^{-1/r} . 
\ee
where   
$C $ depends only on $r$
and $\omega$.

Moreover, if \textsc{Alsm} terminates in (A1), there holds
\be
\label{outRHS2}
\#\t\cT_k-\#\cT_k \le C  c(f)^{ 1/r} \e^{-1/r},  
\ee
with the same dependence of $C$ as above.
\end{lemma}
\noindent
{\bf Proof:} Suppose that  in step (A1) $\delta \le \omega \| f_\delta - \cB \bar u_\delta\|_\Om$,
so that (A1) outputs $f_k, \t\cT_k, \bar u_k$. If this occurs when having entered (A1) for the first time (i.e. $\delta_0=\delta = \|f\|_\Om$), 
we
observe that $\|f-  f_\delta\|_\Om\le \delta$ implies $\|  f_\delta\|_\Om\le 2\delta$, hence 
$$
\|  f_\delta- \cB \bar u_\delta\|_\Om \le \|  f_\delta\|_\Om + \|\bar u_\delta\|_\U\le 3\delta,
$$
by stability. Now
 $r$-optimality says that 
$$
\#\t\cT_0-\#\cT_0 \le 3^{1/r}c(f)^{1/r} (3\delta)^{-1/r}
\le 3^{1/r}c(f)^{1/r} \| f_k- \cB \bar u_k\|_\Om^{-1/r},
$$
which confirms the assertion in this case. If $(A1)$ terminates at the $p$th  step of the loop \eref{while2} for some $p>1$ and let $\t\cT_{k,q}$, $f_{\delta,q}, \bar u_{\delta,q}$   denote the 
quantities, generated at the $q$th stage 
in that loop   for some $q\le p$.
Then, for $q<p$, one has $\omega\| f_{\delta,q}-\cB \bar u_{\delta, q}\|_\Om < 2^{-q}\delta_{k}$
while $\omega\| f_{\delta,q}-\cB \bar u_{\delta, q}\|_\Om \ge 2^{-p}\delta_{k}$. On the other hand, \eref{foptimal}
  implies that
$$
\begin{aligned}
\#\t\cT_{k,p} - \#\cT_k &= \#\t\cT_{k,p} -\#\t\cT_{k,p-1} + \cdots + \#\t\cT_{k, 1}- \#\cT_k \\
&\le c(f)^{1/r}\sum_{q=0}^{p } (2^{-q}\delta_k)^{-1/r} 
 \le  c(f)^{1/r}2^{\frac{p-1}r}\delta_k^{-1/r},
\end{aligned}
$$
as claimed under \eref{outRHS}. 

To confirm \eref{outRHS2}, observe that
$$
\e < 2^{-(p-1)}\delta_k + \|f_{k,p-1}-\cB \bar u_{k,p-1}\|_\Om,
$$
while also $\omega\|f_{k,p-1}-\cB \bar u_{k,p-1}\|_\Om\le 2^{-(p-1)}\delta_k$. Hence
$\e\le (2^{-(p-1)}+\omega^{-1})\delta_k$ whence the assertion follows. \hfill $\Box$\\

The statement of the following lemma is also a standard step in the complexity analysis of \textsc{Afem}.
However, the argument needs to be modified because an error reduction can only be asserted
after a group of subdivision steps. 

\begin{lemma}
\label{lem:basicest}
Given $\theta\in (0,1)$, assume that $(u^*,f)\in \cA^r$. For $\rho<1$ according to \eref{solution}, fixed
let $\zeta,\omega>0$ be  constants such that \eref{omega} holds and  
\be
\label{omega2}
1- \rho^{-2}(\zeta + (2+\zeta)\omega)^2 \ge\theta.
\ee
Given any admissible $\cT=\t\cT_k$, let $f_\cT=f_k\in \PP_\cT$ and $\bar u_\cT=\textsc{Solve}[f_\cT,\cT]\in \U_\cT$ be the (exact)
solution of the discrete problem \eref{normalT} with right hand side $f_\cT$.
If  
\be
\label{squeeze}
 \|f-f_\cT\|_\Om   
\le \delta_k\le \omega \|f_\cT- \cB \bar u_\cT\|_\Om  
\ee
then, the output $\bar u_\cTr$, $\cTr$ of $\textsc{Reduce}[\cT,f_\cT,\bar u_\cT, \theta,\rho]$ satisfies
\be
\label{count}
\#\cTr - \#\cT \le  
 C\|u^*-\bar u_\cT\|_\U^{-1/r} |u|^{1/r}_{\cA^r},
\ee
with $C$ depending only on $\zeta,\omega$, $\mT$, and the maximal number $m$  of applications
of $\textsc{Refine}$ in step (A2) of  $\textsc{Reduce}[\cT,f_\cT,\bar u_\cT, \theta,\rho]$.

\end{lemma}
\noindent
{\bf Proof:} Given $\cT$ as in the assertion, let $\bar\cT$ be the smallest admissible partion $\bar\cT\succ\cT_0$ such that
$u_{\bar\cT}=\textsc{Solve}[f,\bar\cT]$ satisfies
\be
\label{baru}
\|u^*- u_{\bar\cT}\|_\U\le \zeta\|u^*- \bar u_\cT\|_\U,
\ee
for some $\zeta\in (0,1)$ to be fixed later.
Then, by assumption on $(u^*,f)$,
\be
\label{Ar}
\#\bar\cT - \#\cT_0\le \zeta^{-1/r}\|u^*-\bar u_\cT\|^{-1/r}_\U |u^*|^{1/r}_{\cA^r}.
\ee
Now let $\check\cT$ be the smallest common refinement of $\cT$ and $\bar\cT$. Denote:  

- $u_{\check\cT}\in\U_{\check\cT}$   the (exact) solution of \eref{normalT} (with right
hand side $f$),

-   $\bar u_\cT\in\U_\cT$ is the exact discrete solution in $\U_\cT$
with right hand side $f_\cT$,

-   $\bar u_{\check\cT}\in\U_{\check\cT}$ solves
\eref{normalT} in $\U_{\check\cT}$ with right hand side $f_\cT$,

-  $\bar u^*=\cB^{-1}f_\cT\in\U$ denote the exact solution of \eref{normal} with perturbed right\\
\hspace*{0.8cm} hand side $f_\cT$.\\
Finally, let $\check\Om_{\rm ref}$ denote the refinement region obtained when 
passing from $\cT$ to $\check\cT$. Recall from Theorem \ref{thm:pract1}, (b), that under the above assumptions on $\rho$,   the number $m'$ of refinements needed in
\textsc{Reduc} until \eref{while} fails is bounded by  $ m$. Let us denote the corresponding
intermediate discrete solutions with load $f_\cT$ by $\bar u^j\in \U_{\cT^j}$, for some 
$\cT\prec\cT^j\preceq\cTr$, $j=0,\ldots,m'$, $\bar u^0
=\bar u_\cT $, $\bar u^{m'}= \bar u_\cTr$. Hence
\be
\label{hence}
\|f_\cT- \cB \bar u_\cT^j\|_\Om= \|\bar u^* - \bar u_\cT^j\|_\U> \rho\|\bar u^*- \bar u^{j-1}_\cT\|_\U,\quad
j= 1,\ldots,m'-1.
\ee
We show next that for $\zeta$ sufficiently small $\|\bar u^*- \bar u_{\check\cT}\|_\U\le \lambda 
\|\bar u^*-\bar u_\cT\|_\U$ for some $\lambda <1-\theta$ (see \eref{tilde}) in order to apply then Lemma \ref{lem:reduce}.

To that end, we have, on account of the stability of \eref{normalT},  \eref{baru},
and the fact that $\bar\cT\prec\check\cT$,
\be
\label{mix}
\begin{aligned}
\|\bar u^*- \bar u_{\check\cT}\|_\U &= \|u^*- u_{\check\cT} + u_{\check\cT}- 
\bar u_{\check\cT} 
+ \bar u^*- u^*\|_\U\\
&\le \|u^*- u_{\check\cT}\|_\U + \|\cB_{\check\cT}^{-1}(f- f_\cT)\|_\U +\|f_\cT -f\|_\Omega 
\\
&\le 2\|f_\cT -f\|_\Omega + \|u^*- u_{\check\cT}\|_\U\\
&\stackrel{\eref{baru}}{\le} \zeta \|u^*- \bar u_\cT\|_\U + 2\|f_\cT -f\|_\Omega\\
&= \zeta \|\bar u^* -  \bar u_{\cT}      + u^*- \bar u^* \|_\U + 2\|f_\cT -f\|_\Omega\\
&\le \zeta \|\bar u^* -  \bar u_{\cT} \|_\U +(2+\zeta)\|f_\cT -f\|_\Omega\\
&\stackrel{\eref{squeeze}}{\le}\zeta \|\bar u^* -  \bar u_{\cT} \|_\U+
(2+\zeta)\omega \|\cB \bar u_\cT- f_\cT\|_\Om \\
&= \zeta \|\bar u^* -  \bar u_{\cT} \|_\U + (2+\zeta)\omega \|\bar u^* -  \bar u_{\cT} \|_\U\\
&= (\zeta + (2+\zeta)\omega)\|\bar u^* -  \bar u_{\cT} \|_\U\\
&\stackrel{\eref{hence}}{\le}\rho^{-1}(\zeta + (2+\zeta)\omega)\|\bar u^* -  \bar u^j \|_\U,\quad j<m'.
\end{aligned}
\ee
Choosing $\omega, \zeta$ small enough such that
\be
\label{zeta}
\lambda:=\rho^{-1}(\zeta + (2+\zeta)\omega)^2  
 \le 1-\theta,
\ee
 we have verified \eref{tilde}.  
Denoting by $\check\Om_{\rm ref}$  the refinement region corresponding to $\check\cT\succ \cT$, we apply Lemma   \ref{lem:reduce}  to $f$ replaced by $f_\cT$ to conclude that
\be
\label{bulkcheck}
\|\cB \bar u^j  - f_\cT\|^2_{\check\Om_{\rm ref}}\ge \theta \|\cB \bar u_\cT - f_\cT\|^2_\Om,\quad
j=0,\ldots,m'-1.
\ee 
The number of new cells generated when refining $\cT^{j-1}$ to $\cT^j$, $j=1,\ldots,m'$ is near-minimal.
Hence, there exists a constant $C$ such that
\be
\begin{aligned}
\label{refbound}
\#\cTr - \#\cT^j&\le  C(\#\check\cT- \#\cT_0)\le C(\#\bar\cT- \#\cT_0)\\
&\le C\zeta^{-1/r}\|u^*-\bar u_\cT\|_\U^{-1/r} |u|^{1/r}_{\cA^r}.
\end{aligned}
\ee
Therefore,
\begin{align*}
\#\cTr - \#\cT&\le \#\cTr - \#\cT^{m'-1} +\cdots + \#\cT^1 - \#\cT\\
&\le \sum_{j=0}^{m'-1} \#\cTr - \#\cT^j\\
&\le m  C\zeta^{-1/r}\|u^*-\bar u_\cT\|_\U^{-1/r} |u|^{1/r}_{\cA^r},
\end{align*}
confirming \eref{count}.\hfill $\Box$\\

To complete the proof of Theorem \ref{thm:pract2}, it remains to bound the partition growth. When $f\in \cA^r$ and $(f,\text{Appr})$ is $r$-optimal
(see Assumption {\bf Af}), Lemma \ref{lem:RHS}, \eref{outRHS}, asserts that, at the $k$th stage, the output $[\cT_{k },
f_{k },\bar u_{k }]$ of step 
(A1) satisfies 
\be
\label{A1out}
\#\t\cT_k-\#\cT_k \le C  c(f)^{ 1/r} \delta_k^{-1/r}
\ee
prior to termination. The same inequality holds with $\delta_k$ replaced by $\e$ when the scheme 
terminates in (A1). The constant $C$ depends on $\mT, \zeta, \omega, \rho, m$, satisfying \eref{solution}, \eref{omega}, and \eref{omega2}, noting that $\zeta$ and $\omega$ are chosen in dependence of $\theta$
which, in turn determines $\xi$ in \eref{xi} and thereby $\rho$ and $m$.

We know from the  proof of Theorem \ref{thm:pract1} that  $\textsc{Reduce}[\cT_{k },f_k,\bar u_k,\theta,\rho]$
produces an approximate solution $u_{k+1}\in \U_{\cT_{k+1}}$ such that $\| f_k 
- \cB u_{k+1}\|_\Om\le \rho\| f_k- \cB \bar u_k\|_\Om$, and that
$\bar u_\cT=\bar u_{k}$,   satisfy \eref{squeeze}, as output of step (A1). By Lemma \ref{lem:basicest}, 
$$
\begin{aligned}
\#\cT_{k+1}- \#\t\cT_{k }&\le  
  C|u^*|_{\cA^r}^{1/r}\|u^*- \bar u_{k}\|^{-1/r}_\U, 
\end{aligned}
$$
where $C$ depends on $\zeta,\mT,\omega,m$.
Since 
\begin{align*}
\|u^*- \bar u_k\|_\U&= \|u^*- u^*_k+ u^*_k-\bar u_k\|_\U=\|f-f_k+ f_k- \cB \bar u_k\|_\Om\\
&\ge \|f_k- \cB \bar u_k\|_\Om- \|f- f_k\|_\Om,
\end{align*}
and recalling that when exiting step (A1), we have $\|f- f_k\|_\Om\le \delta_k$ and
$\delta_k\le \omega\|f_k- \cB \bar u_k\|_\Om$,
we conclude that $\|u^*- \bar u_k\|_\U\ge (\omega^{-1}-1)\delta_k$ and therefore
\be
\label{deltak}
\#\cT_{k+1}- \#\t\cT_{k }\le  C (\omega^{-1}-1)^{-1/r} |u^*|_{\cA^r}^{1/r}\delta_k^{-1/r}.
\ee
Combining this with \eref{A1out}, yields
\be
\label{Tcombine}
\begin{aligned}
\#\cT_{k+1}-\#\cT_k&\le  
  C'\big(c(f)^{1/r}+ |u|_{\cA^r}^{1/r}\Big)
\delta_k^{-1/r} .
\end{aligned}
\ee
where $C'$ depends on $\zeta,\omega,\mT, \rho,\theta,m$.\hfill $\Box$

%%%%%%%%%%%%%%%%%%%%%%%%%%%%%%%%%%%%
\section{Concluding Remarks}
%%%%%%%%%%%%%%%%%%%%%%%%%%%%%%%

As in the case of classical \textsc{Afem}, D\"orfler marking appears to   be the right way
to ensure linear convergence and, under natural assumptions, also optimal rates in terms of
partition growth. Yet there are a few noteworthy distinctions between the convergence analyses 
in both scenarios. While the appearance of local meshsizes in the a posteriori error indicators,
used in classical \textsc{Afem}, allows one to establish a fixed error reduction per adaptive cycle 
\eref{cycle}, it seems that in \textsc{Alsm} this can only be shown under very stringent assumptions
on a certain new variant of data-visibility (see \eref{reliable}). The key is therefore to better understand the interrelation between this data-visibility term and the current residual which
can be done with the aid of projections to several auxiliary subspaces of $\U$, related to the
current local refinement. This allows one to establish a fixed error reduction  for a finite number
of refinement steps - in brief ``block reduction'', provided that the load $f$ is sufficiently well
aligned with the image under $\cB$ of certain complement spaces, associated with the refinements.
Simply {\em assuming} this to hold (like assuming that data-oscillation is of higher order) would
give a relatively simple error reduction result (see \eref{reliable} in Theorem \ref{thm:main1}).
The deeper insight lies therefore in tying ``visibility'' of $f$ to its piecewise polyomial approximability
which is a variant of a regularity property.
The resulting linear
convergence in terms   of a ``block reduction'' remains the key to a   complexity analysis, 
which, however, requires a somewhat amended reasoning concerning complexity. It is perhaps interesting
that then there don't seem to be any constraints on the bulk parameter $\theta$ 
other than being less than one. 

Natural future goals  are to explore to what extend the above concepts carry over to minimum residual
methods when the residual is measured in a non-trivial dual space instead of $\L$, or 
in yet a different direction to {\em nonlinear} FOSLS formulations.

%%%%%%%%%%%%%%%%%%%%%%%%%%%%%%%%

%%%%%%%%%%%%%%%%%%%%%%%

\begin{thebibliography}{99}


\bibitem{BDO}
M. Bachmayr, W. Dahmen, M. Oster,   {Variationally correct neural residual regression for parametric {PDE}s: On the viability of controlled accuracy},
 {IMA Journal of Numerical Analysis}, (2025),
  doi{10.1093/imanum/draf073}
 
\bibitem{BD}
P. Binev, R. DeVore, Fast computation in adaptive tree approximation, Numer. Math., 97 (2004),  193--217.

\bibitem{BDD}
 P. Binev, W. Dahmen, R. DeVore, Adaptive Finite Element Methods with Convergence Rates,  Numer. Math., 97(2004), 219--268.


\bibitem{FOSLS}
 Z. Cai,   R. Lazarov,  T. A. Manteuffel,   S. F. McCormick,
 {First-Order System Least Squares for Second-Order Partial Differential Equations: Part {I}},
{SIAM Journal on Numerical Analysis}, 31 (1994),
 {1785--1799}.
 
\bibitem{5authors}
W. Dahmen,   I. Graham, B. Faehrmann, W. Hackbusch, S. Sauter, 
Inverse Inequalities on Non-Quasiuniform Meshes
and Application to the Mortar Element Method,  
 Math. Comp. 73  (2004), 1107-1138. 
 

\bibitem{FK} 
T. F\"{u}hrer,   M. Karkulik,
 {Space-time least-squares finite elements for parabolic equations},
 {Comput. Math. Appl.},
  {92} (2021), 27--36.

 

\bibitem{FGK}
 T. F\"{u}hrer,R. Gonzales, M. Karkulik,   
  {Well-Posedness of First Order Acoustic Wave Equations and Space-Time Finite Element Approximation},
 {arXiv preprint arXiv:2311.10536},
   {2023}
 

\bibitem{CDG}
 P. Cort\'es Castillo,W. Dahmen, J Gopalakrishnan, 
  {DPG loss functions for learning parameter-to-solution maps by neural networks},
 {arXiv preprint: http://arxiv.org/abs/2506.18773},
  {2025}
 

\bibitem{Bochev}
Pavel Bochev,  and Max D. Gunzburger,  
 {Least-squares finite element methods},
 {Applied Mathematical Sciences},
 {166},
 {Springer, New York},
 {2009},
ISBN  {978-0-387-30888-3},
 {https://doi.org/10.1007/b13382}.
 
\bibitem{BCNVacta}
A. Bonito, C. Canuto, R. Nochetto, A Veeser, Acta Numerica (2024),   163--485,    doi:10.1017/S0962492924000011
 
\bibitem{BPS}
D. Braess, V. Pillwein,  J. Sch\"oberl, Equilibrated residual error estimates are
p-robust, Comput. Methods Appl. Mech. Engrg 198  (2009), 1189--1197.

\bibitem{Braess}
D. Braess, Finite Elements: Theory, Fast Solvers, and Applications in Solid Mechanics,
Cambridge University Press,
  2010,
ISBN:
9780511618635,
doi.org/10.1017/CBO9780511618635

\bibitem{BBF}
D. Boffi, F. Brezzi, M. Fortin, Mixed Finite Element Methods and Applications, Springer Series in Computational Mathematics, vol. 44, Springer, 2013.

\bibitem{Bringmann}
P. Bringmann, Review and computational comparison of adaptive least-squares finite element schemes, Computers \& Mathematics with Applications, 122 (2024), 1--15.

\bibitem{CP}
C. Carstensen,   E.-J. Park. Convergence and optimality of adaptive least squares finite element methods. SIAM J. Numer. Anal., 53 (2015), 43--62.

\bibitem{axioms}
C. Carstensen, M. Feischl, M. Page,   D. Praetorius, Axioms of adaptivity, Comput. Math. Appl., 67  (2014), 1195--1253.

\bibitem{DGacta}
L. Demkowicz, J. Gopalakrishnan,
The discontinuous Petrov-Galerkin method, Acta Numerica (2025), pp. 293-384,
 doi:10.1017/S0962492924000102


\bibitem{Feischl}
M. Feischl, Optimality of a standard adaptive finite element method for the Stokes problem, SIAM J. Numer. Anal., 57  (2019), 1124--1157.

\bibitem{GS}
G. Gantner, R. Stevenson,
Further results on a space-time FOSLS formulation of parabolic PDEs,
ESAIM: Mathematical Modelling and Numerical Analysis, 55 (2021), 283--299.

\bibitem{KV}
C. Kreuzer,  A. Veeser, Oscillation in a posteriori error estimation, Numer. Math. 148  (2021), 43--78.

\bibitem{QCD}
Y. Qiu, P. Chen, W. Dahmen,
Variationally correct operator learning: Reduced basis neural operator
  with a posteriori error estimation, 
   Computer Methods in Applied Mechanics and Engineering, 461 (2026), 119228,
   doi: 10.1016/j.cma.2026.119228,
   Preprint Dec 2025, arXiv:2512.21319, [math.NA], [cs.LG], [cs.NA]

\bibitem{robfocm}
R. Stevenson, Optimality of a standard adaptive finite element method, Found. Comput. Math., 7 (2007), 
  245--269.
  
\bibitem{S}
R. Stevenson, The completion of locally refined simplicial partitions created by bisection, Math. Comp.,
 77  (2008), 227--241.

\bibitem{MNS}
P. Morin, R. H. Nochetto, K. G. Siebert, Convergence of adaptive finite element methods, SIAM Rev. 44,  (2002) 631--658. 


\bibitem{NSV}
R. H. Nochetto, K. G. Siebert,   A. Veeser, Theory of adaptive finite element methods: An introduction, in Multiscale, Nonlinear and Adaptive Approximation (R. DeVore
and A. Kunoth, eds), Springer,  2009, pp. 409--542.

\bibitem{V} R. Verf\"{u}rth, {\it A Review of A Posteriori Error Estimation and Adaptive
Mesh-Refinement Techniques},
 Wiley-Teubner, Chichester, 1996.

\end{thebibliography}
\end{document}